\documentclass[12pt]{amsart}

\usepackage{amsmath,amsfonts,amssymb,graphicx,hyperref, enumitem,amsthm,mathabx,amsrefs}
\usepackage{array}
\hypersetup{colorlinks=true,citecolor=blue,linkcolor=blue,urlcolor=blue,pdfstartview=FitH}

\newcommand{\paren}[1]{\ensuremath{\left( #1 \right)}}
\newcommand{\set}[1]{\ensuremath{\left\{ #1 \right\}}}
\newcommand{\innerprod}[1]{\ensuremath{\left< #1 \right>}}
\newcommand{\abs}[1]{\ensuremath{\left| #1 \right|}}
\newcommand{\setdiv}{\,\middle|\,}

\newcommand{\Z}{\mathbb{Z}}
\newcommand{\R}{\mathbb{R}}
\newcommand{\N}{\mathbb{N}}
\newcommand{\Q}{\mathbb{Q}}
\newcommand{\C}{\mathbb{C}}

\newcommand{\wbar}[1]{\overline{#1}}
\newcommand{\wtilde}[1]{\widetilde{#1}}
\newcommand{\what}[1]{\widehat{#1}}
\newcommand{\wcheck}[1]{\widecheck{#1}} 
\newcommand{\e}[1]{e\paren{#1}}
\newcommand{\pFqName}[2]{{_{#1}F_{#2}}}
\newcommand{\pFq}[5]{\pFqName{#1}{#2}\paren{\begin{matrix}#3;\\#4;\end{matrix}\,#5}}
\newcommand{\trans}[1]{{#1}^T}

\newcommand{\Matrix}[1]{\begin{pmatrix}#1\end{pmatrix}}

\newcommand{\Min}[1]{\ensuremath{\min \set{#1}}}

\DeclareMathOperator*{\res}{res}
\DeclareMathOperator{\diag}{diag}

\theoremstyle{plain} 
\newtheorem{thm}{Theorem}
\newtheorem{cor}[thm]{Corollary}
\newtheorem{lem}[thm]{Lemma}
\newtheorem{prop}[thm]{Proposition}

\newcommand{\piecewise}[1]{\left\{\begin{matrix}#1\end{matrix}\right.}
\newcommand{\If}{\mbox{if }}

\renewcommand{\Re}{{\mathop{\mathgroup\symoperators Re}}}
\renewcommand{\Im}{{\mathop{\mathgroup\symoperators Im}}}
\newcommand{\sgn}{{\mathop{\mathgroup\symoperators \,sgn}}}

\DeclareMathAlphabet{\mathcalligra}{T1}{calligra}{m}{n}
\newcommand{\WigDName}{\mathcal{D}}
\newcommand{\WigDMat}[1]{\WigDName^{#1}}

\newcommand{\WigD}[3]{\WigDMat{#1}_{#2,#3}}
\newcommand{\WigdName}{\mathcalligra{d}}
\newcommand{\Wigd}[3]{\WigdName^{#1}_{#2,#3}}

\newcommand{\Dtildek}[2]{\mathcal{R}^{#1}\paren{#2}}

\newcommand{\AdSq}{{\mathrm{Ad}}^2}

\newcommand{\vpmpm}[1]{v_{_{#1}}}

\newcommand{\bv}{\frak{v}}

\newcommand{\specmu}{\mathbf{spec}}
\newcommand{\cosmu}{\mathbf{cos}}
\newcommand{\sinmu}{\mathbf{sin}}

\title{Kuznetsov, Petersson and Weyl on GL(3), II:\\ The generalized principal series forms.}
\author{Jack Buttcane}
\date{26 August 2019}
\address{Department of Mathematics \& Statistics, 5752 Neville Hall, Orono, ME 04469, USA}
\email{jack.buttcane@maine.edu}
\thanks{During the time of this research, the author was supported by NSF grant DMS-1601919.}

\begin{document}

\begin{abstract}
This paper initiates the study by analytic methods of the generalized principal series Maass forms on $GL(3)$.
These forms occur as an infinite sequence of one-parameter families in the two-parameter spectrum of $GL(3)$ Maass forms, analogous to the relationship between the holomorphic modular forms and the spherical Maass cusp forms on $GL(2)$.
We develop a Kuznetsov trace formula attached to these forms at each weight and use it to prove an arithmetically-weighted Weyl law, demonstrating the existence of forms which are not self-dual.
Previously, the only such forms that were known to exist were the self-dual forms arising from symmetric-squares of $GL(2)$ forms.
The Kuznetsov formula developed here should take the place of the $GL(2)$ Petersson trace formula for theorems ``in the weight aspect''.
As before, the construction involves evaluating the Archimedian local zeta integral for the Rankin-Selberg convolution and proving a form of Kontorovich-Lebedev inversion.
\end{abstract}


\maketitle

\section{Introduction}

The non-principal, generalized principal series forms for $GL(3)$ are forms of minimal $K$-type attached to the $(2d+1)$-dimensional Wigner-$\WigDName$ matrix $\WigDMat{d}$ (see section \ref{sect:Notation}) with $d \ge 2$ and spectral parameters
\[ \mu = \mu(r) := \paren{\tfrac{d-1}{2}+r,-\tfrac{d-1}{2}+r,-2r}. \]
These generate strict subrepresentations of principal series representations which are induced from representations on the $2,1$ Levi subgroup.
As a one-parameter family of Maass forms, one might compare them to the spherical Maass cusp forms on $GL(2)$, but as a lower-dimensional subspace of the full $GL(3)$ spectrum, one might also compare them to the point spectrum, i.e. holomorphic modular forms, on $GL(2)$.
The purpose of this paper is to initiate the analytic aspects of their study.
This is the last of the three spectral Kuznetsov formulae on $GL(3)$ (for full level over $\Q$), the others were constructed in \cite{SpectralKuz} and \cite{WeylI}.
From the adelic perspective, we are handling the case of ramification at the place at infinity and this would be largely unaffected by considering quotients by congruence subgroups, aka. ramification at the finite places.

The structure of this paper is similar to the previous one, with the chain of constructions
\[ \text{Stade's formula }\Rightarrow \text{ Kontorovich-Lebedev inversion } \Rightarrow \text{ Kuznetsov's formula } \Rightarrow \text{ Weyl law}, \]
but the technical details of producing an infinite sequence of formulae are much more involved.
In particular, the proof of Stade's formula (Theorem \ref{thm:Stade} below) becomes difficult, and this is the piece we were unable to complete in the preceeding paper; the reader will notice the construction of section \ref{sect:StadesFormula} is quite intricate.
On the other hand, the proof of the Weyl law becomes even easier since the Mellin-Barnes integrals for the Kuznetsov kernel functions are much simplified over the principal series case.

Let $\mathcal{S}^{d*}_3$ be a basis of vector-valued minimal-weight cusp forms attached to the $(2d+1)$-dimensional representation of $SO(3,\R)$.
Denote the spectral parameter of such a Maass form $\varphi$ by $r_\varphi \in \C$ and notice the Ramanujan-Selberg conjecture $r_\varphi \in i\R$ is known for these forms \cite[Theorem 3]{HWII}.
\begin{thm}
\label{thm:HeckeWeyl}
	Let $0 \ne r'\in i\R$ and $d \ge 3$, $T > M  > 1$ such that $d+T\to\infty$.
	Then if we assume (as we may) that $\mathcal{S}^{d*}_3$ consists of Hecke eigenforms,
	\begin{align*}
		\sum_{\abs{r_\varphi-T r'} < M} \frac{1}{L(1,\AdSq \varphi)} =& \frac{3}{2\pi} \int_{\abs{r-T r'} < M} \specmu^d(r) dr+O\paren{d(d+T)^2},
	\end{align*}
	where $\specmu^d(r)$ is the spectral weight
	\begin{align*}
		\specmu^d(r) =& \tfrac{1}{16\pi^4 i}(d-1)\paren{\tfrac{d-1}{2}-3r}\paren{\tfrac{d-1}{2}+3r}.
	\end{align*}
\end{thm}
The main term in the Weyl law is $\asymp dM(d+T)^2$, and we fail to give an asymtotic exactly when $M \ll 1$, but a more careful analysis in that range would tighten up the error term to save a small power of $\log(d+T)$ over the main term (so the asymptotic holds for $M$ down to a slightly negative power of $\log(d+T)$), as the argument is actually using a Gaussian of width $(\log(d+T))^{-1/2}$; we leave this to future explorers.
One can see from \eqref{eq:TestFunOffSymmBound} that this is the limit of the method (up to $\log \log$ factors, etc.), and this (roughly) conforms to the rule of thumb that the Kuznetsov formula cannot resolve a ball of radius less than one in the spectrum.

From the Weyl law, we see that the generalized principal series forms are fewer in number than the principal series forms of the previous papers:
If we identify the weight $d=0$ and weight $d=1$ tempered spectra with $\set{\mu\in\C^3\setdiv \Re(\mu)=0,\mu_1+\mu_2+\mu_3=0}$, then those Weyl laws (see \cite{WeylI}; the $d=0$ case follows almost verbatim using \cite{SpectralKuz} as it appears in \cite[Theorem 10.0]{ArithKuzI}) may be written
\[ \sum_{\substack{\varphi\in \mathcal{S}^{d*}_3\\ \mu_\varphi \in T\Omega}} \frac{1}{L(1,\AdSq \varphi)} = \frac{3}{2\pi} \int_{T\Omega} \specmu^d(\mu) d\mu+O\paren{T^{4+\epsilon}}, \]
where $\Omega$ is some nice subset of the spectrum, $d\mu=d\mu_1 \, d\mu_2$, and the spectral weights are
\begin{align*}
	\specmu^0(\mu) =& \frac{1}{384\pi^4} \prod_{i<j} (\mu_i-\mu_j)\tan\frac{\pi}{2}(\mu_i-\mu_j), \\
	\specmu^1(\mu) =& \frac{1}{64\pi^4} \paren{\prod_{i<j} (\mu_i-\mu_j)} \cot\frac{\pi}{2}(\mu_1-\mu_3) \cot\frac{\pi}{2}(\mu_2-\mu_3) \tan\frac{\pi}{2}(\mu_1-\mu_2).
\end{align*}
The main terms of the two principal series Weyl laws are $\asymp T^5$, while an equivalent statement for fixed $d \ge 3$ would have main term $\asymp T^3$.
However, in every case, the size of the spectral measure is (generically) $\asymp \prod_{i<j} \abs{\mu_i-\mu_j}$, and a sum of generalized principal series forms over $d \asymp T$ restores equality between the two types.
This should not be too surprising, as the number of holomorphic modular forms of weight $k \asymp T$ is also of similar size to the number of spherical $GL(2)$ Maass cusp forms of spectral parameter $\abs{\mu} \asymp T$.

Prior to the above theorem, existence theorems for Maass forms on $GL(3)$ were limited to two types:
First, M\"uller's Weyl law \cite[Theorem 0.1]{Muller01} indicates the existence of $GL(3)$ Maass forms by counting all forms having a given $K$-type (including lifts of lower-weight forms).
Second, the symmetric-square lift construction of Gelbart and Jacquet \cite{GelbartJacquet} directly produces a $GL(3)$ Maass form from a given $GL(2)$ Maass form (usually in the adelic viewpoint, so computing the effect of ramification in the classical sense can be trying).
The group of forms studied here is too small to be detected by M\"uller's theorem (they are drowned out by the lifts of the $d=0,1$ forms; see previous paragraph), and the symmetric square lifts of $GL(2)$ holomorphic modular forms of weight $k$ are known to occur at $d=2k-1$ with $r=0$, so for $\varphi \in \mathcal{S}^{d*}_3$ when the weight $d \ge 3$ is even or when the spectral parameter $r_\varphi$ is non-zero, the existence of such forms is new.

The unweighted Weyl laws (i.e. without $1/L(1,\AdSq \varphi)$) can be obtained from the Kuznetsov formulae given below by a weight-removal technique, see \cite{MeFan01}.
One then expects the spectral measures $\specmu^d(\mu) d\mu$ to correspond to the measure in Harish-Chandra's Plancherel theorem for $SL(3,\R)$ (see \cite[Theorem 13.4.1]{Wallach}), up to some predictable (but difficult) constant, and this is easy to see for the spherical case, as the measure is nicely written out in \cite[Theorem 4.3.1(1)]{T02}; meanwhile, the remaining cases $d \ge 1$ are somewhat difficult to extract from the literature. 

The theorem above also does not handle the case $d=2$ for technical reasons, and we will discuss this in section \ref{sect:Weight2}, where we provide a rough upper bound on the Weyl law for those forms (which is needed in \cite{ArithKuzII}).

\section{Some notation for $GL(3)$}
\label{sect:Notation}
Let $G=PSL(3,\R) = GL(3,\R)/\R^\times$ and $\Gamma=SL(3,\Z)$.
The Iwasawa decomposition of $G$ is $G=U(\R) Y^+ K$ using the groups $K=SO(3,\R)$,
\[ U(R) = \set{\Matrix{1&x_2&x_3\\&1&x_1\\&&1} \setdiv x_i\in R}, \qquad R \in \set{\R,\Q,\Z}, \]
\[ Y^+ = \set{\diag\set{y_1 y_2,y_1,1} \setdiv y_1,y_2 > 0}. \]
The measure on the space $U(\R)$ is simply $dx := dx_1 \, dx_2 \, dx_3$, and the measure on $Y^+$ is
\[ dy := \frac{dy_1 \, dy_2}{(y_1 y_2)^3}, \]
so that the measure on $G$ is $dg :=dx \, dy \, dk$, where $dk$ is the Haar probability measure on $K$ (see \cite[sect. 2.2.1]{HWI}).
We generally identify elements of quotient spaces with their coset representatives, and in particular, we view $U(\R)$, $Y^+$, $K$ and $\Gamma$ as subsets of $G$.

Characters of $U(\R)$ are given by
\[ \psi_m(x) = \psi_{m_1,m_2}(x) = \e{m_1 x_1+m_2 x_2}, \qquad \e{t} = e^{2\pi i t}, \]
where $m\in\R^2$; we say $\psi=\psi_m$ is non-degenerate when $m_1 m_2 \ne 0$.
Characters of $Y^+$ are given by the power function on $3\times 3$ diagonal matrices, defined by
\[ p_\mu\paren{\diag\set{a_1,a_2,a_3}} = \abs{a_1}^{\mu_1} \abs{a_2}^{\mu_2} \abs{a_3}^{\mu_3}, \]
where $\mu\in\C^3$.
We assume $\mu_1+\mu_2+\mu_3=0$ so this is defined modulo $\R^\times$, renormalize by $\rho=(1,0,-1)$, and extend by the Iwasawa decomposition
\[ p_{\rho+\mu}\paren{x y k} = y_1^{1-\mu_3} y_2^{1+\mu_1}, \qquad x\in U(\R),y\in Y^+,k\in K. \]
Integrals in $\mu$ use the permutation-invariant measure $d\mu=d\mu_1 \, d\mu_2$.

The Weyl group $W$ of $G$ contains the six matrices
\begin{equation*}
	\begin{array}{rclcrclcrcl}
		I &=& \Matrix{1\\&1\\&&1}, && w_2 &=& -\Matrix{&1\\1\\&&1}, && w_3 &=& -\Matrix{1\\&&1\\&1}, \\
		w_4 &=& \Matrix{&1\\&&1\\1}, && w_5 &=& \Matrix{&&1\\1\\&1}, && w_l &=& -\Matrix{&&1\\&1\\1},
	\end{array}
\end{equation*}
with the relations $w_3 w_2=w_4$, $w_2 w_3=w_5$ and $w_2 w_3 w_2=w_3 w_2 w_3=w_l$.
The Weyl group induces an action on the coordinates of $\mu$ by $p_{\mu^w}(a):=p_\mu(waw^{-1})$, and we denote the coordinates of the permuted parameters by $\mu^w_i := (\mu^w)_i$, $i=1,2,3$.
Explicity, the action is
\begin{align*}
	\mu^I =& \paren{\mu_1,\mu_2,\mu_3}, & \mu^{w_2} =& \paren{\mu_2,\mu_1,\mu_3}, & \mu^{w_3} =& \paren{\mu_1,\mu_3,\mu_2}, \\
	\mu^{w_4} =& \paren{\mu_3,\mu_1,\mu_2}, & \mu^{w_5} =& \paren{\mu_2,\mu_3,\mu_1}, & \mu^{w_l} =& \paren{\mu_3,\mu_2,\mu_1}.
\end{align*}

The group of diagonal, orthogonal matrices $V \subset G$ contains the four matrices $v_{\varepsilon_1,\varepsilon_2}=\diag\set{\varepsilon_1,\varepsilon_1 \varepsilon_2,\varepsilon_2}, \varepsilon\in\set{\pm 1}^2$, which we abbreviate $V = \set{\vpmpm{++}, \vpmpm{+-}, \vpmpm{-+}, \vpmpm{--}}$.
We write $Y=Y^+ V$ for the diagonal matrices in $G$ so that $Y \cap K = V$.
We tend not to distinguish between elements of $Y$ and pairs in $(\R^\times)^2$, as the multiplication is the same.

We will also require the Bruhat decomposition $G=U(\R) Y W U(\R)$; the decomposition becomes unique if we replace, for each element $w\in W$, the right-hand copy of $U(\R)$ with $\wbar{U}_w(\R)$ where $\wbar{U}_w = (w^{-1} \trans{U} w) \cap U$.
When taking the Bruhat decomposition of an element $\gamma\in\Gamma$, we have $\gamma=bcvwb'$ with $b,b'\in U(\Q)$, $v\in V$, $w\in W$ and $c$ of the form
\begin{align*}
	\Matrix{\frac{1}{c_2}\\&\frac{c_2}{c_1}\\&&c_1}, \qquad c_1, c_2\in \N.
\end{align*}
For such a matrix and a pair of characters $\psi_m,\psi_n$ $m,n\in\Z^2$, we define the Kloosterman sum attached to the Weyl element $w\in W$ by
\[ S_w(\psi_m,\psi_n,c) = \sum_{\substack{b cw b' \in U(\Z)\backslash \Gamma/(V \wbar{U}_w(\Z))}} \psi_m(b) \psi_n(b'), \]
provided it is well-defined, and zero otherwise; these are given explicitly in \cite[sect. 4.2]{WeylI}.

The unitary, irreducible representations of $K$, up to isomorphism, are given by the Wigner $\WigDName$-matrices $\WigDMat{d}:K\to GL(2d+1,\C)$ for $0 \le d \in \Z$.
We treat $\WigDMat{d}$ primarily as a matrix-valued function with the usual properties
\[ \WigDMat{d}(kk')=\WigDMat{d}(k)\WigDMat{d}(k'), \qquad \trans{\wbar{\WigDMat{d}(k)}} = \WigDMat{d}(k)^{-1} = \WigDMat{d}(k^{-1}). \]
The entries of the matrix-valued function $\WigDMat{d}$ are indexed from the center:
\[ \WigDMat{d} = \Matrix{\WigD{d}{-d}{-d}&\ldots&\WigD{d}{-d}{d}\\ \vdots&\ddots&\vdots\\ \WigD{d}{d}{-d}&\ldots&\WigD{d}{d}{d}}. \]
The entries, rows, and columns of of the derived matrix- and vector-valued functions (e.g. the matrix-valued Whittaker function, see section \ref{sect:WhittFuns}) will be indexed similarly.
As the Wigner $\WigDName$-matrices exhaust the equivalence classes of unitary, irreducible representations of the compact group $K$, they give a basis of $L^2(K)$, as in \cite[sect. 2.2.1]{HWI}, by the Peter-Weyl theorem.
We tend to refer to the index $d$ as the ``weight'' of the Wigner $\WigDName$-matrix and any associated objects (e.g. the Whittaker function, Maass forms, etc.).
The Wigner $\WigDName$-matrices on the $v$-matrices occur frequently, and these are given by (see \cite[section 2.2.2]{HWI})
\begin{align*}
	\WigDMat{d}(v_{\varepsilon,+1}) =& \diag\set{\varepsilon^d,\ldots,\varepsilon^{-d}}, & \WigD{d}{m'}{m}(v_{\varepsilon,-1}) =& (-1)^d\varepsilon^{m'} \delta_{m'=-m}.
\end{align*}

A Maass form (cuspidal or Eisenstein) of weight $d$ and spectral parameters $\mu$ for $\Gamma$ is a row vector-valued (or matrix-valued) smooth function $f:\Gamma\backslash G\to\C^{2d+1}$ that transforms under $K$ as $f(gk)=f(g)\WigDMat{d}(k)$, satisfies a moderate growth condition and is an eigenfunction of both Casimir operators with eigenvalues matching $p_{\rho+\mu}$ (see \cite[sect. 3]{HWII}).
We denote the quadratic Casimir operator by $\Delta_1$ and the cubic Casimir operator by $\Delta_2$; precise expressions for these operators will not be necessary, here.
The action of the Lie algebra of $G$ on Maass forms gives rise to five operators $Y^a$, $\abs{a} \le 2$ on vector-valued Maass forms (see \cite[sect. 5.2]{HWII}) that change the weight $d \mapsto d+a$, and a Maass form is said to have minimal weight if it is sent to zero by the lowering operators $Y^{-1}, Y^{-2}$ and an eigenfunction of the $Y^0$ operator.

Throughout the paper, we take the term ``smooth'', in reference to some function, to mean infinitely differentiable on the domain.
The letters $x,y,k,g, v$ and $w$ will generally refer to elements of $U(\R)$, $Y^+$, $K$, $G$, $V$, and $W$, respectively.
The letter $\psi$ will generally refer to a character of $U(\R)$, and $\mu$ will always refer to an element of $\C^3$ satisfying $\mu_1+\mu_2+\mu_3=0$.
Vectors (resp. matrices) not directly associated with the Wigner $\WigDName$-matrices, e.g. elements $n \in \Z^2$ are indexed in the traditional manner from the left-most entry (resp. the top-left entry), e.g.  $n=(n_1,n_2)$.
We do not use the primed notation $F'$ for derivatives, but rather to distinguish functions and variables with similar purpose.

\section{Results}
Throughout we assume $d \ge 2$ is an integer.
Suppose $r\in i\R$ and $y =\diag(y_1 y_2,y_1,1)$ in $Y^+ \cong (\R^+)^2$, the space of positive diagonal matrices as in \cite[sect. 2.1]{HWI}.
Then for $\abs{m'}\le d$, we write $m'=\varepsilon m$ with $\varepsilon = \pm 1$ and $0 \le m \le d$, and define
\begin{align}
\label{eq:MinWhittGPS}
	W^{d*}_{\varepsilon m}(y, r) =& \frac{1}{2^{d+1}\pi} \sqrt{\binom{2d}{d+m}} \sum_{\ell=0}^{m} \varepsilon^\ell \binom{m}{\ell} \int_{\Re(s)=\mathfrak{s}} (2\pi y_1)^{1-s_1} (2\pi y_2)^{1-s_2} \\
	& \times \Gamma\paren{\tfrac{d-1}{2}+s_1-r} \Gamma\paren{\tfrac{d-1}{2} + s_2+r} B\paren{\tfrac{d-m+s_1+2r}{2},\tfrac{\ell+s_2-2r}{2}} \frac{ds}{(2\pi i)^2}, \nonumber
\end{align}
for any $\mathfrak{s} \in (\R^+)^2$.
The completed Whittaker function $W^{d*}(y,r)$ attached to the minimal $K$-type $\WigDMat{d}$ is the row vector with coordinates $\paren{W^{d*}_{-d},\ldots,W^{d*}_d}$.

We start in section \ref{sect:StadesFormula} by generalizing Stade's formula to the Rankin-Selberg convolution of two generalized principal series forms of the same minimal weight (see \cite[Corollary 3]{WeylI}).
\begin{thm}
\label{thm:Stade}
	Define
	\[ \Psi^d=\Psi^d(r,r',t) = \int_{Y^+} W^{d*}(y,r) \trans{W^{d*}(y,r')} (y_1^2 y_2)^t dy, \]
	then
	\begin{align*}
		\Psi^d=& 2^{4-d-4t-r-r'} \pi^{2-3t} \Gamma\paren{t+r+r'} \Gamma\paren{\tfrac{d-1}{2}+t+r-2r'} \Gamma\paren{\tfrac{d-1}{2}+t+r'-2r} \\
		& \times \Gamma(d-1+t+r+r') \Gamma\paren{\tfrac{t}{2}-r-r'} / \Gamma\paren{\tfrac{3t}{2}}.
	\end{align*}
\end{thm}
In the preceeding paper, we essentially followed Stade's original proof \cite{Stade02}, but that fails here for technical reasons; the proof of this result is dramatically more complicated, and this is the central, new step that we were unable to accomplish previously.

Define the spectral weights
\begin{align*}
	\frac{1}{\sinmu^d(r)} :=& \frac{2\pi i}{3} \lim_{t\to 0^+} t \, \Psi^d(r,-r,t) = 2^{5-d} \pi^3 i \Gamma(d-1) \Gamma\paren{\tfrac{d-1}{2}-3r} \Gamma\paren{\tfrac{d-1}{2}+3r}, \\
	\frac{1}{\cosmu^d(r)} :=& \Psi^d(r,-r,1)= \frac{2^{1-d}}{\pi} \Gamma(d) \Gamma\paren{\tfrac{d+1}{2}+3r} \Gamma\paren{\tfrac{d+1}{2}-3r}, \\
	\specmu^d(r) :=& \frac{\sinmu^d(r)}{\cosmu^d(r)} = \tfrac{1}{16\pi^4 i}(d-1)\paren{\tfrac{d-1}{2}-3r}\paren{\tfrac{d-1}{2}+3r}.
\end{align*}

In section \ref{sect:KontLeb}, we apply Stade's formula to generalize Kontorovich-Lebedev inversion to our choice of Whittaker functions, using the method of Goldfeld and Kontorovich \cite{GoldKont}.
\begin{thm}
\label{thm:KontLebedev}
	For $f:Y^+\to\C^{2d+1}$ define
	\[ f^\sharp(r) = \int_{Y^+} f(y) \wbar{\trans{W^{d*}(y,r)}} dy, \]
	and for $F:i\R\to\C$, define
	\[ F^\flat(y) = \int_{\Re(r)=0} F(r) W^{d*}(y,r) \sinmu^d(r)\,dr. \]
	If $F(r)$ is holomorphic and Schwartz-class on a vertical strip $\set{r\setdiv\abs{\Re(r)}<\delta}$ for some $\delta > 0$, then
	\[ (F^\flat)^\sharp(r) = F(r). \]
\end{thm}
Once again, we are making no claim as to the image of $F \mapsto F^\flat$ beyond the necessary convergence.

In section \ref{sect:KuzForm}, the Kuznetsov formula will follow from Stade's formula and Kontorovich-Lebedev inversion.
The Kuznetsov kernel functions are defined from the power-series solutions \eqref{eq:JwlDef},  \eqref{eq:Jw4Def} as the linear combinations
\begin{align}
\label{eq:KwIEval}
	 K_I^d(y;r) =& 1,
\end{align}
\begin{align}
\label{eq:Kw4Eval}
	 &4\pi \cos\pi\paren{\tfrac{d}{2}+3r} K_{w_4}^d(y;r) =\\
	 & \qquad (-\varepsilon_1 i)^d J_{w_4}(y,\mu(r)^{w_4}) \exp \paren{-\tfrac{\varepsilon i\pi}{2}\paren{\tfrac{d-1}{2}-3r}} -(-\varepsilon_1 i)^d J_{w_4}(y,\mu(r)), \nonumber
\end{align}
\begin{align}
\label{eq:Kw5}
	K_{w_5}^d(y;r) =& K_{w_4}^d(\vpmpm{-+} y^\iota;-r),
\end{align}
\begin{align}
\label{eq:KwlEval}
	-4\pi \cos\pi\paren{\tfrac{d}{2}+3r} K_{w_l}^d(y,r) =& \delta_{\varepsilon_1=-1} (\varepsilon_2)^d J_{w_l}(y,\mu(r)^{w_4}) + \delta_{\varepsilon_2=-1} (-\varepsilon_1)^d J_{w_l}(y,\mu(r)) \\
	&\qquad -\delta_{\varepsilon_1 \varepsilon_2=-1} (-\varepsilon_1)^d J_{w_l}(y,\mu(r)^{w_3}). \nonumber
\end{align}
These functions are the kernels for the integral transforms
\begin{align}
\label{eq:HwDef}
	H_w(F; y) =& \frac{1}{\abs{y_1 y_2}} \int_{\Re(r)=0} F(r) K^d_w(y, r) \specmu^d(r) dr.
\end{align}

Let $\mathcal{S}_3^{d*}$ be an orthonormal basis of $SL(3,\Z)$ cusp forms of minimal weight $d$ and $\mathcal{S}_2^{d*}$ an orthonormal basis of $SL(2,\Z)$ cusp forms of minimal weight $d$ (i.e. those coming from holomorphic modular forms of weight $d$).
Note that $\mathcal{S}_2^{d*}$ is empty for odd $d$ and the symmetric-squares of holomorphic modular forms \textit{only} occur for odd $d$, but we will see there are forms in $\mathcal{S}_3^{d*}$ with $r \ne 0$ for all $d$.
The particular normalization of the Fourier-Whittaker coefficients $\rho_\varphi^*(m)$ of $\varphi\in\mathcal{S}_3^{d*}$ and $\rho_\phi^*(m;r)$ of the maximal parabolic Eisenstein series attached to $\phi\in\mathcal{S}_2^{d*}$ is given in \eqref{eq:FWcoefDef}.

We are now ready to state the Kuznetsov formula for this class of forms.
\begin{thm}
\label{thm:Kuznetsov}
	Let $d \ge 3$.
	Let $F(r)$ be Schwartz-class and holomorphic on $\abs{\Re(r)}<\frac{1}{4}+\delta$ for some $\delta > 0$.
	Then for $m,n\in\Z^2$ with $m_1 m_2 n_1 n_2 \ne 0$,
	\[ \mathcal{C}+\mathcal{E}=\mathcal{K}_I+\mathcal{K}_4+\mathcal{K}_5+\mathcal{K}_l, \]
	where
	\begin{align*}
		\mathcal{C} =& \sum_{\varphi \in \mathcal{S}_3^{d*}} F(r_\varphi) \frac{\wbar{\rho_\varphi^*(m)} \rho_\varphi^*(n)}{\cosmu^d(r_\varphi)}, \\
		\mathcal{E} =&\frac{2}{2\pi i} \sum_{\phi \in \mathcal{S}_2^{d*}} \int_{\Re(r)=0} F(r) \frac{\wbar{\rho_\phi^*(m;r)} \rho_\phi^*(n;r)}{\cosmu^d(r)} dr,
	\end{align*}
	\begingroup
	\allowdisplaybreaks
	\begin{align*}
		\mathcal{K}_I =& \delta_{\substack{\abs{m_1}=\abs{n_1}\\ \abs{m_2}=\abs{n_2}}} H_I(F;(1, 1)) \nonumber \\
		\mathcal{K}_4 =& \sum_{\varepsilon\in\set{\pm1}^2} \sum_{\substack{c_1,c_2\in\N\\ \varepsilon_1 m_2 c_1=n_1 c_2^2}} \frac{S_{w_4}(\psi_m,\psi_{\varepsilon n},c)}{c_1 c_2} H_{w_4}\paren{F; \paren{\varepsilon_1 \varepsilon_2 \tfrac{m_1 m_2^2n_2}{c_2^3n_1},1}}, \\
		\mathcal{K}_5 =& \sum_{\varepsilon\in\set{\pm1}^2} \sum_{\substack{c_1,c_2\in\N\\ \varepsilon_2 m_1 c_2= n_2 c_1^2}} \frac{S_{w_5}(\psi_m,\psi_{\varepsilon n},c)}{c_1 c_2} H_{w_5}\paren{F; \paren{1,\varepsilon_1 \varepsilon_2 \tfrac{m_1^2 m_2 n_1}{c_1^3 n_2}}}, \\
		\mathcal{K}_l =& \sum_{\varepsilon\in\set{\pm1}^2} \sum_{c_1,c_2\in\N} \frac{S_{w_l}(\psi_m,\psi_{\varepsilon n},c)}{c_1 c_2} H_{w_l}\paren{F; \paren{\varepsilon_2 \tfrac{m_1 n_2 c_2}{c_1^2}, \varepsilon_1 \tfrac{m_2 n_1 c_1}{c_2^2}}}.
	\end{align*}
	\endgroup
\end{thm}
The absolute convergence of the Kuznetsov formula can be shown through contour shifting and relies on the holomorphy condition; we will demonstrate this explicitly in the proof of the Weyl law.

It will be useful to compare the current kernel functions to those for the other minimal $K$-types.
Using the $K^{\pm\pm}_{w_l}$ functions of \cite[sect. 4.6]{WeylI} and \cite{SpectralKuz}, and setting $\varepsilon=\sgn(y)$, we have
\begin{align}
	-4\pi \cos\pi\paren{\tfrac{d}{2}+3r} K_{w_l}^d(y,r) =\piecewise{
		0 & \If \varepsilon=(1,1), \\
		K^{-+}_{w_l}(y,\mu(r)^{w_4}) & \If \varepsilon=(-1,1),\\
		(-1)^d K^{+-}_{w_l}(y,\mu(r)) & \If \varepsilon=(1,-1),\\
		K^{--}_{w_l}(y,\mu(r)) & \If \varepsilon=(-1,-1).}
\end{align}
On the other hand, we can provide a much more compact statement of the Mellin-Barnes integrals at $\mu=\mu(r)$; define
\begin{align*}
\begin{array}{*3{>{\displaystyle}l}}
	B^{--}_{w_l}(s,r) :=& (-1)^d B\paren{s_1+3r,s_2-3r}, \\
	B^{-+}_{w_l}(s,r) :=& B\paren{s_2-3r,1-s_1-s_2}, \\
	B^{+-}_{w_l}(s,r) :=& B\paren{s_1+3r,1-s_1-s_2}, \\
	B^{++}_{w_l}(s,r) :=& 0,
\end{array}&&
	Q(d,s) :=& \frac{\Gamma\paren{\frac{d-1}{2}+s}}{\Gamma\paren{\frac{d+1}{2}-s}},
\end{align*}
then
\begin{align}
\label{eq:KwlMB}
	K_{w_l}^d(y,r) =& \frac{1}{4\pi^2} \abs{\frac{y_2}{y_1}}^r \int_{-i\infty}^{+i\infty} \int_{-i\infty}^{+i\infty} \abs{4\pi^2 y_1}^{1-s_1} \abs{4\pi^2 y_2}^{1-s_2} B^{\varepsilon}_{w_l}\paren{s,r} Q(d,s_1) Q(d,s_2) \frac{ds_2}{2\pi i} \frac{ds_1}{2\pi i}.
\end{align}
And of course, the $K_{w_4}^d$ integral is new once again:
\begin{align}
\label{eq:Kw4MB}
	 K_{w_4}^d(y;\mu(r)) =& \frac{4^{r-1}(\varepsilon_1 i)^d}{\pi^{\frac{3}{2}}} \int_{-i\infty}^{+i\infty} \abs{4\pi^3 y_1}^{1-r-s} Q(d,s) \paren{\frac{\Gamma\paren{\frac{s+3r}{2}}}{\Gamma\paren{\frac{1-s-3r}{2}}}+\varepsilon_1 i \frac{\Gamma\paren{\frac{1+s+3r}{2}}}{\Gamma\paren{\frac{2-s-3r}{2}}}} \frac{ds}{2\pi i} \\
	 =& \frac{(\varepsilon_1 i)^d }{4\pi^2} \int_{-i\infty}^{+i\infty} \abs{8\pi^3 y_1}^{1-r-s} Q(d,s) \Gamma\paren{s+3r} \exp\paren{\varepsilon_1 \tfrac{i\pi}{2}\paren{s+3r}} \frac{ds}{2\pi i}. \nonumber
\end{align}
The author would like to point out that this is essentially the only possibility for a Mellin-Barnes integral (on a vertical contour) for any linear combination of the $J_{w_4}$ (which is a $\pFqName02$ hypergeometric function), which lends some credence to the hope that we have not missed any $v$-matrices along the way.

In section \ref{sect:TechnicalWeylLaw}, we prove our first application of the Kuznetsov formula, which is a technical version of the Weyl law with both analytic and arithmetic weights:
\begin{thm}
\label{thm:TechnicalWeyl}
	Let $F(r)$ and $\delta$ be as in Theorem \ref{thm:Kuznetsov}, then we have
	\begin{align*}
		\sum_{\varphi \in \mathcal{S}^{d*}_3} F(r_\varphi) \frac{\abs{\rho_\varphi^*(1)}^2}{\cosmu^d(r_\varphi)} = \int_{\Re(r)=0} F(r) \specmu^d(r) dr+O\paren{E_1+E_2},
	\end{align*}
	where
	\begin{align*}
		E_1 :=& \int_{\Re(r)=0} \abs{F(r)} (d+\abs{r})^{1+\epsilon} \abs{dr}, \\
		E_2 :=& \int_{\Re(r)=-\frac{1}{4}-\eta} \paren{\abs{F(r)}+\abs{F(-r)}} d (d+\abs{r})^{-\frac{1}{2}+\epsilon} \abs{dr},
	\end{align*}
	and we assume $0 < \eta < \delta$ satisfies $\eta = O(\epsilon)$ and $\epsilon>0$ arbitrarily small.
\end{thm}
In section \ref{sect:WeylLaw}, by a careful choice of test function, we will remove the analytic weights:
\begin{cor}
\label{cor:Weyl}
	For $\Omega$,  $r'$ and $T$ and $M$ as in Theorem \ref{thm:HeckeWeyl},
	\[ \sum_{\abs{r_\varphi-T r'} < M} \frac{\abs{\rho_\varphi^*(1)}^2}{\cosmu^d(r_\varphi)} = \int_{\abs{r-T r'} < M} \specmu^d(r) dr+O\paren{d(d+T)^2}. \]
\end{cor}
We have made no attempt to optimize the error terms in Theorem \ref{thm:TechnicalWeyl}, as this will be sufficient for what is likely the best-possible error term in the corollary due to the sharp cut-off.

A Rankin-Selberg argument using Stade's formula gives the Kuznetsov formula on Hecke eigenvalues, and this is discussed in section \ref{sect:RS}.
\begin{thm}
\label{thm:HeckeKuznetsov}
	When the bases of Theorem \ref{thm:Kuznetsov} are taken to be Hecke eigenfunctions, the left-hand side may be written as
	\begin{align*}
		\mathcal{C} =& \frac{2\pi}{3} \sum_{\varphi \in \mathcal{S}_3^{d*}} F(r_\varphi) \frac{\wbar{\lambda_\varphi(m)} \lambda_\varphi(n)}{L(1,\AdSq \varphi)}, \\
		\mathcal{E} =&\frac{2}{i} \sum_{\phi \in \mathcal{S}_2^{d*}} \int_{\Re(r)=0} F(r) \frac{\wbar{\lambda_\phi(m,r)} \lambda_\phi(n,r)}{L(\phi,1+3r) L(\phi,1-3r) L(1,\AdSq \phi)} dr.
	\end{align*}
	where $\lambda_\varphi(m)$ and $\lambda_\phi(m,r)$ as in \cite[(5.13)]{HWI} are the Hecke eigenvalues of the associated forms.
\end{thm}

The main theorem now follows.

\section{Background}
\subsection{The Whittaker functions}
\label{sect:WhittFuns}
The Whittaker function $W^{d*}(y, r)$ naturally extends to a function on $G$ by the Iwasawa decomposition (see \cite[sect. 2.4]{HWI}),
\[ W^{d*}(xyk, r) = \psi_{1,1}(x) W^{d*}(y, r) \WigDMat{d}(k). \]
In \cite[Theorem 6]{HWII}, this Mellin-Barnes integral derives from the matrix-valued Whittaker function
\begin{align}
\label{eq:JacWhittDef}
	W^d(g,\mu,\psi) := \int_{U(\R)} I^d(w_l u g, \mu) \wbar{\psi(u)} du, \qquad I^d(xyk,\mu) := p_{\rho+\mu}(y) \WigDMat{d}(k),
\end{align}
in the form
\begin{align}
\label{eq:WhittCompletion}
	\Lambda^*(r) W^d_{-d}(y, \mu(r), \psi_{1,1}) =& W^{d*}(y, r), \qquad W^d_{d,m} = 0,
\end{align}
where
\begin{align*}
	\Lambda^*(r) =\Lambda^*(\mu(r)) =& \frac{(-1)^d}{\pi} (2\pi)^{-\tfrac{d-1}{2}-3r} \Gamma(d) \Gamma\paren{\tfrac{d+1}{2}+3r}.
\end{align*}

The functional equations in $\mu$ of the matrix-valued Whittaker function are given by \cite[prop. 3.3]{HWI},
\begin{align}
\label{eq:WhittFEs}
	W^d(g,\mu,\psi_{1,1}) = T^d(w,\mu) W^d(g,\mu^w,\psi_{1,1}), \qquad w \in W,
\end{align}
and these are generated by the matrices
\begin{align}
\label{eq:Tdw2}
	T^d(w_2,\mu) :=&\pi^{\mu_1-\mu_2} \Gamma^d_\mathcal{W}(\mu_2-\mu_1,+1), \\
\label{eq:Tdw3}
	T^d(w_3,\mu) :=& \pi^{\mu_2-\mu_3} \WigDMat{d}(\vpmpm{--}w_l) \Gamma^d_\mathcal{W}(\mu_3-\mu_2,+1) \WigDMat{d}(w_l\vpmpm{--}),
\end{align}
where $\Gamma^d_\mathcal{W}(u,\varepsilon)$ is a diagonal matrix coming from the functional equation of the classical Whittaker function \cite[(2.20)]{HWI}:
If $\mathcal{W}^d(y,u)$ is the diagonal matrix-valued function with entries (see \cite[sect. 2.3.1]{HWI})
\begin{equation}
\label{eq:classWhittDef}
\begin{aligned}
	\mathcal{W}^d_{m,m}(y,u) =& \int_{-\infty}^\infty \paren{1+x^2}^{-\frac{1+u}{2}} \paren{\frac{1+ix}{\sqrt{1+x^2}}}^{-m} \e{-yx} dx \\
	=& \piecewise{\displaystyle \frac{(\pi\abs{y})^{\frac{1+u}{2}}}{\abs{y} \Gamma\paren{\frac{1-\varepsilon m+u}{2}}} W_{-\frac{\varepsilon m}{2}, \frac{u}{2}}(4\pi\abs{y}) & \If y \ne 0, \\[5pt]
\displaystyle \frac{2^{1-u}\pi \,\Gamma(u)}{\Gamma\paren{\frac{1+u+m}{2}} \Gamma\paren{\frac{1+u-m}{2}}} & \If y=0,}
\end{aligned}
\end{equation}
(where $W_{\alpha,\beta}(y)$ is the classical Whittaker function), then for $y\ne 0$, we have the functional equations
\begin{align}
\label{eq:WhittGammas}
	\mathcal{W}^d(y,-u) =& (\pi\abs{y})^{-u} \Gamma^d_\mathcal{W}(u,\sgn(y)) \mathcal{W}^d(y,u) & \Gamma_{\mathcal{W},m,m}^d(u,\varepsilon) =& \frac{\Gamma\paren{\frac{1-\varepsilon m+u}{2}}}{\Gamma\paren{\frac{1-\varepsilon m-u}{2}}}.
\end{align}

The matrices $T^d(w,\mu)$ commute with the action of the $v$-matrices according to \cite[(3.26)]{HWI}
\begin{align}
\label{eq:WhittFEvsV}
	\WigDMat{d}(v) T^d(w,\mu) = T^d(w,\mu) \WigDMat{d}(w^{-1} v w).
\end{align}

We can express $\Gamma_\mathcal{W}^d(u,+1)$ in terms of $\mathcal{W}^d(0,-u)$ and visa-versa by \cite[(130)]{HWII} and its inverse
\begin{align}
\label{eq:GammaWtoW0}
	\Gamma_\mathcal{W}^d(u,+1) =& \frac{i\,\Gamma(1+u)}{2^{1+u}\pi} \paren{\exp\paren{\tfrac{i\pi u}{2}} \WigDMat{d}(\vpmpm{--})-\exp\paren{-\tfrac{i\pi u}{2}} \WigDMat{d}(\vpmpm{+-})} \WigDMat{d}(w_2)\mathcal{W}^d(0,-u), \\
\label{eq:W0toGammaW}
	\mathcal{W}^d(0,u) =& 2^{-u} \Gamma(u) \paren{\exp\paren{\tfrac{i\pi u}{2}} \WigDMat{d}(\vpmpm{--})+\exp\paren{-\tfrac{i\pi u}{2}} \WigDMat{d}(\vpmpm{+-})} \WigDMat{d}(w_2) \Gamma_\mathcal{W}^d(-u,+1).
\end{align}

\subsection{The Spectral Expansion}
There exists a symmetric differential operator $\Lambda_{\frac{d-1}{2}}$ of \cite[prop. 2]{HWII}, whose kernel in the weight $d$ forms is exactly the span of the forms with spectral parameters of the form $\mu(r)$, $r\in\C$.
Then for a Schwartz-class function $f:\Gamma\backslash G\to\C^{2d+1}$ satisfying
\begin{align}
\label{eq:TestFunReqs}
	f(gk) =& f(g)\WigDMat{d}(k), & \Lambda_{\frac{d-1}{2}} f = 0,
\end{align}
the spectral expansion of \cite{HWI, HWII} takes the form
\begin{align}
\label{eq:SpectralExpand}
	f(g) =& \sum_{\varphi \in \mathcal{S}_3^{d*}} \varphi(g) \, \int_{\Gamma\backslash G} f(g') \wbar{\trans{\varphi(g')}} dg' \\
		& \qquad +\frac{2}{2\pi i} \sum_{\phi \in \mathcal{S}_2^{d*}} \int_{\Re(r)=0} E^d_d(g, \phi, r) \, \int_{\Gamma\backslash G} f(g') \wbar{\trans{E^d_d(g', \phi, r)}} dg' \, dr, \nonumber
\end{align}
where $\mathcal{S}_3^{d*}$ is a basis of $SL(3,\Z)$ cusp forms of minimal weight $d$ and $\mathcal{S}_2^{d*}$ is a basis of $SL(2,\Z)$ cusp forms of minimal weight $d$.
Note that there are no $SL(2,\Z)$ cusp forms and hence no Eisenstein series when $d$ is odd.
We have dropped all terms whose spectral parameters are not of the form $\mu(r)$ by the usual orthogonality argument:
\[ 0=\innerprod{\Lambda_{\frac{d-1}{2}} f, g} = \innerprod{f, \Lambda_{\frac{d-1}{2}} g} = \lambda_g\paren{\Lambda_{\frac{d-1}{2}}} \innerprod{f,g}, \]
whenever $g$ is an eigenfunction of the Casimir operators.
The extra 2 in the coefficient of the maximal parabolic Eisenstein series comes from the equality
\[ E^d_d(g', \phi, r) = E^d_{-d}(g', \phi, r). \]
The cusp forms $\varphi\in\mathcal{S}_3^{d*}$ are normalized by $\left<\varphi,\varphi\right>=1$, in place of \cite[(166)]{HWII}.

\subsection{Integrals of gamma functions and hypergeometric identities}
\label{sect:GammaInts}
We require a number of identities of Mellin-Barnes integrals, and we collect them here for ease of disposition.
Let $\mathcal{C}$ denote any contour from $-i\infty$ to $i\infty$ which obeys the Barnes integral convention (no gamma function in the numerator should have its argument pass through the negative real axis) and allow $\mathcal{C}$ to vary from line to line.
The parameters $a,b,c,d\in\C$ and arguments $x,y\in\R$, $z\in\C$ may be any values such that corresponding integrals and hypergeometric series converge absolutely and avoid the relevant branch cuts.
We will not need to worry about the branch cuts, but more precise statements can be found in the references.

The hypergeometric functions in general have the integral description \cite[7.2.3.12-13]{Prud3}
\begin{align}
\label{eq:pFqToMB}
	\frac{\prod_{i=1}^p \Gamma(a_i)}{\prod_{i=1}^q \Gamma(b_i)} \pFq{p}{q}{a_1,\ldots,a_p}{b_1,\ldots,b_q}{-z} =& \int_{\mathcal{C}} \frac{\Gamma(s)\prod_{i=1}^p \Gamma(a_i-s)}{\prod_{i=1}^q \Gamma(b_i-s)} z^{-s} \frac{ds}{2\pi i},
\end{align}
for $q \le p+1$.

By Mellin inversion and the definition of the beta function,
\begin{align}
\label{eq:BetaInvMellin}
	\int_{\mathcal{C}} \frac{\Gamma(s) \Gamma(a+1)}{\Gamma(a+1+s)} x^{-s} \frac{ds}{2\pi i} &= \piecewise{(1-x)^a & \If 0<x<1, \\ 0 & \If x \ge 1.}
\end{align}

From \eqref{eq:pFqToMB} and Thomae's theorem \cite[7.4.4.2]{Prud3}
\begin{align}
\label{eq:ThomaeMB}
	& \int_{\mathcal{C}} \frac{\Gamma(s) \Gamma\paren{a-s} \Gamma\paren{a+b-d-s}}{\Gamma\paren{2-c+s} \Gamma\paren{a+b-s}\Gamma\paren{a+c-d-s}} \frac{ds}{2\pi i} \\
	&= \frac{\Gamma(a)}{\Gamma(2-c) \Gamma(d) \Gamma\paren{c-b}} \int_{\mathcal{C}} \frac{\Gamma(s) \Gamma(1-s) \Gamma(d-s) \Gamma\paren{c-b-s}}{\Gamma\paren{c-s} \Gamma\paren{1+a-s}} \frac{ds}{2\pi i}. \nonumber
\end{align}

The Euler integral representations
\begin{align}
\label{eq:EulerBeta}
	B(u,v) =& \int_0^1 x^{u-1} (1-x)^{v-1} dx,
\end{align}
and \cite[7.2.1.2]{Prud3}
\begin{align}
\label{eq:Elem2F1}
	\int_0^1 (1-x)^a x^b (1+y x)^c dx =& \frac{\Gamma(a+1)\Gamma(b+1)}{\Gamma(a+b+2)} \pFq21{b+1,-c}{a+b+2}{-y}.
\end{align}

The Pfaff transformation \cite[7.2.1.7]{Prud3} at $z=\frac{1}{2}$
\begin{align}
\label{eq:Pfaff}
	\pFq21{a,b}{c}{\tfrac{1}{2}} &= 2^{-b} \pFq21{a,c-b}{c}{-1}.
\end{align}

Barnes' first lemma \cite[7.2.1.4]{Prud3}
\begin{align}
\label{eq:BarnesFirst}
	\int_{\mathcal{C}} \Gamma(a+s)\Gamma(b+s)\Gamma(c-s)\Gamma(d-s) \frac{ds}{2\pi i} = \frac{\Gamma(a+c)\Gamma(a+d)\Gamma(b+c)\Gamma(b+d)}{\Gamma(a+b+c+d)}.
\end{align}

Gauss' theorem \cite[7.3.5.2]{Prud3}
\begin{align}
\label{eq:GaussThm}
	\pFq21{a,b}{c}{1}=\frac{\Gamma(c)\Gamma(c-a-b)}{\Gamma(c-a)\Gamma(c-b)}.
\end{align}

And several more, apparently unnamed, identities:
\begin{align}
\label{eq:SimpleMB}
	\int_{\mathcal{C}} x^{-s/2} \Gamma\paren{a+s} \Gamma\paren{b-s} \frac{ds}{2\pi i} =& \paren{1+\sqrt{x}}^{-a-b} x^{a/2} \, \Gamma(a+b), && \text{\cite[8.4.2.5]{Prud3}},
\end{align}
\begin{align}
\label{eq:MB1F0Eval}
	 \int_{\mathcal{C}} \frac{\Gamma(s-a) \Gamma\paren{b-s}}{\Gamma\paren{b-a}} \frac{ds}{2\pi i} =& \pFq10{b-a}{}{-1} = 2^{a-b}, && \text{\cite[7.3.1.1]{Prud3}},
\end{align}
\begin{align}
\label{eq:2F1toBeta}
	\frac{1}{a} \pFq21{a+b,a}{a+1}{-1}+\frac{1}{b} \pFq21{a+b,b}{b+1}{-1} =& B(a,b), && \text{\cite[7.3.5.5]{Prud3}}.
\end{align}

\section{Stade's Formula}
\label{sect:StadesFormula}
The proof of Theorem \ref{thm:Stade} is quite complicated and involves a number of seemingly random manipulations for which the author has no intuition beyond their simple effectiveness.
To assist the reader, we give a brief summary:
The proof proceeds by the usual application of Parseval's formula for the Mellin transform in \eqref{eq:MinWhittGPS} followed by expanding the beta functions using Euler's integral \eqref{eq:EulerBeta}, at which point we can evaluate all of the sums and inverse Mellin transforms to produce a single, elementary, two-dimensional integral \eqref{eq:PsidElementary} for $\Psi^d$.
The elementary transformation \eqref{eq:ElementarySubs} splits the integral into two, simpler integrals and Mellin-expanding the resulting hypergeometric integrals gives $\Psi^d$ as a sum of two, three-dimensional Mellin-Barnes integrals at 1 in \eqref{eq:PsidSimpleMB}.
Applying Thomae's theorem allows us to evaluate one of the three one-dimensional integrals, at which point the hypergeometric identities \eqref{eq:Pfaff} and \eqref{eq:2F1toBeta} recombine the sum of two, two-dimensional integrals into one, one-dimensional integral, which can be evaluated by Barnes' first lemma.

It seems that the hypergeometric manipulations would be better realized as substitutions on the elementary integral; in particular, splitting the integral into two pieces only to later recombine them suggests we missed an elementary substitution which kept the two pieces together in the first place.
Of course, even knowing such a thing should exist doesn't necessarily make it easy to find.

We denote the contour for an $n$-dimensional Mellin-Barnes integral as $\mathcal{C}^n$, continuing on from section \ref{sect:GammaInts} for simplicity of notation.

Now to the proof:
Starting from the definition of $\Psi^d$, we apply the definition of the completed Whittaker function \eqref{eq:MinWhittGPS}, and Parseval's formula for the Mellin transform so that
\begin{align*}
	\Psi^d=& \frac{(2\pi)^{4-3t}}{2^{2d+2}\pi^2} \int_{\mathcal{C}^2} F(s) \Gamma\paren{\tfrac{d-1}{2}+s_1-r} \Gamma\paren{\tfrac{d-1}{2} + s_2+r} \\
	& \times \Gamma\paren{\tfrac{d-1}{2}+2t-s_1-r'} \Gamma\paren{\tfrac{d-1}{2} + t-s_2+r'} \frac{ds}{(2\pi i)^2}, \\
	F(s) :=& \sum_{m'=-d}^d \binom{2d}{d+m'} \sum_{\ell_1=0}^{m} \varepsilon^{\ell_1} \binom{m}{\ell_1} \sum_{\ell_2=0}^{m} \varepsilon^{\ell_2} \binom{m}{\ell_2} B\paren{\tfrac{d-m+s_1+2r}{2},\tfrac{\ell_1+s_2-2r}{2}} \\
	& \times B\paren{\tfrac{d-m+2t-s_1+2r'}{2}, \tfrac{\ell_2+t-s_2-2r'}{2}},
\end{align*}
after collecting the sums of beta functions.
Now we apply the Euler integral \eqref{eq:EulerBeta} for each beta function, and evaluate the sums, using the binomial theorem in the forms
\begin{align*}
	\sum_{\ell=0}^{m} \varepsilon^{\ell} \binom{m}{\ell} (1-x)^{\frac{\ell}{2}} =& \paren{1+\varepsilon \sqrt{1-x}}^m, \\
	\sum_{m'=-d}^d \binom{2d}{d+m'} u^{m'} =& \paren{2+u^{-1}+u}^d.
\end{align*}
And so $F$ is given by the two-dimensional integral
\begin{align*}
	F(s) =& 2^d \int_0^1 \int_0^1 x_1^{\frac{s_1+2r}{2}-1} (1-x_1)^{\frac{s_2-2r}{2}-1} x_2^{\frac{2t-s_1+2r'}{2}-1}(1-x_2)^{\frac{t-s_2-2r'}{2}-1} \\
	& \times \paren{1+\sqrt{x_1x_2}+\sqrt{1-x_1}\sqrt{1-x_2}}^d dx_1 \, dx_2,
\end{align*}
where we have used several times the fact that
\[ \frac{\sqrt{x}}{1+\sqrt{1-x}}=\frac{1-\sqrt{1-x}}{\sqrt{x}}. \]

Returning to $\Psi^d$, we have
\begin{align*}
	\Psi^d=& \frac{(2\pi)^{4-3t}}{2^{d+2}\pi^2} \int_{[0,1]^2} x_1^{r-1} (1-x_1)^{-r-1} x_2^{t+r'-1} (1-x_2)^{\frac{t-2r'}{2}-1} \paren{1+\sqrt{x_1x_2}+\sqrt{1-x_1}\sqrt{1-x_2}}^d \\
	& \times \int_{\mathcal{C}} \paren{\frac{x_2}{x_1}}^{-s_1/2} \Gamma\paren{\tfrac{d-1}{2}+s_1-r} \Gamma\paren{\tfrac{d-1}{2}+2t-s_1-r'} \frac{ds_1}{2\pi i} \\
	& \times \int_{\mathcal{C}} \paren{\frac{1-x_2}{1-x_1}}^{-s_2/2} \Gamma\paren{\tfrac{d-1}{2} + s_2+r} \Gamma\paren{\tfrac{d-1}{2} + t-s_2+r'} \frac{ds_2}{2\pi i} dx.
\end{align*}

Next we apply \eqref{eq:SimpleMB} to achieve the elementary integral description,
\begin{align}
\label{eq:PsidElementary}
	\Psi^d=& \frac{(2\pi)^{4-3t}}{2^{d+2}\pi^2} \Gamma(d-1+2t-r-r') \Gamma(d-1+t+r+r') \\
	& \times \int_{[0,1]^2} x_1^{\frac{\frac{d-1}{2}+2t+2r-r'}{2}-1} (1-x_1)^{\frac{\frac{d-1}{2} +t+r'-2r}{2}-1} x_2^{\frac{\frac{d-1}{2}+2t+2r'-r}{2}-1} (1-x_2)^{\frac{\frac{d-1}{2} +t+r-2r'}{2}-1} \nonumber \\
	& \times \paren{1+\sqrt{x_1x_2}+\sqrt{1-x_1}\sqrt{1-x_2}}^d \nonumber \\
	& \times \paren{\sqrt{x_1}+\sqrt{x_2}}^{-(d-1+2t-r-r')} \paren{\sqrt{1-x_1}+\sqrt{1-x_2}}^{-(d-1+t+r+r')} dx. \nonumber
\end{align}

Now split the integral at $x_1 = x_2$, and perform the substitutions
\begin{align}
\label{eq:ElementarySubs}
	\left\{\begin{array}{lllll}
		x_1 \mapsto x_1^2 x_2, &\text{then}& \dfrac{1-x_2}{1-x_1^2 x_2} \mapsto x_2^2&\text{on}& x_1<x_2,\\[10pt]
		x_2 \mapsto x_1 x_2^2, &\text{then}& \dfrac{1-x_1}{1-x_1 x_2^2} \mapsto x_1^2&\text{on}& x_2<x_1.
	\end{array}\right.
\end{align}
We have
\begin{align*}
	\Psi^d=& \frac{(2\pi)^{4-3t}}{2^d \pi^2} \Gamma(d-1+2t-r-r') \Gamma(d-1+t+r+r') \\
	& \times \Biggl(\int_{[0,1]^2} x_1^{\frac{d-3}{2}+2t+2r-r'} (1-x_1)^{\frac{t}{2}-r-r'-1} (1+x_1)^{-\frac{3t}{2}} x_2^{\frac{d-3}{2}+t+r-2r'} \\
	& \times (1-x_2)^{t+r+r'-1} (1-x_1 x_2)^{1-\frac{3t}{2}} (1+x_1 x_2)^{1-d-\frac{3t}{2}} dx \\
	& + \int_{[0,1]^2} x_1^{\frac{d-3}{2}+t+r'-2r} (1-x_1)^{t+r+r'-1} x_2^{\frac{d-3}{2}+2t+2r'-r} \\
	& \times (1-x_2)^{\frac{t}{2}-r-r'-1} (1+x_2)^{-\frac{3t}{2}} (1-x_1 x_2)^{1-\frac{3t}{2}} (1+x_1 x_2)^{1-d-\frac{3t}{2}} dx\Biggr).
\end{align*}
Notice the three-summand $d$-th power in \eqref{eq:PsidElementary} has factored due to the substitutions.

Temporarily assuming $0<t<\frac{2}{3}$, we apply \eqref{eq:BetaInvMellin} to the factor $(1-x_1 x_2)^{1-\frac{3t}{2}}$, so that we may apply \eqref{eq:Elem2F1} and \eqref{eq:pFqToMB} twice to produce
\begin{align}
\label{eq:PsidSimpleMB}
	\Psi^d=& \frac{(2\pi)^{4-3t}}{2^d \pi^2} \frac{\Gamma(d-1+2t-r-r') \Gamma(d-1+t+r+r') \Gamma\paren{2-\frac{3t}{2}} \Gamma\paren{\frac{t}{2}-r-r'}}{\Gamma\paren{d-1+\frac{3t}{2}} \Gamma\paren{\frac{3t}{2}}} \\
	& \times \Gamma\paren{t+r+r'} \int_{\mathcal{C}^3} \frac{\Gamma(s_1) \Gamma(s_2) \Gamma(s_3) \Gamma\paren{d-1+\frac{3t}{2}-s_2} \Gamma\paren{\frac{3t}{2}-s_3}}{\Gamma\paren{2-\frac{3t}{2}+s_1}} \nonumber \\
	& \times \Biggl(\frac{\Gamma\paren{\frac{d-1}{2}+t+r-2r'-s_1-s_2} \Gamma\paren{\frac{d-1}{2}+2t+2r-r'-s_1-s_2-s_3}}{\Gamma\paren{\frac{d-1}{2}+2t+2r-r'-s_1-s_2} \Gamma\paren{\frac{d-1+5t}{2}+r-2r'-s_1-s_2-s_3}} \nonumber \\
	& + \frac{\Gamma\paren{\frac{d-1}{2}+t+r'-2r-s_1-s_2} \Gamma\paren{\frac{d-1}{2}+2t+2r'-r-s_1-s_2-s_3}}{\Gamma\paren{\frac{d-1}{2}+2t+2r'-r-s_1-s_2} \Gamma\paren{\frac{d-1+5t}{2}+r'-2r-s_1-s_2-s_3}}\Biggr) \frac{ds}{(2\pi i)^3} \nonumber
\end{align}

Now applying Thomae's theorem in the form \eqref{eq:ThomaeMB} to the $s_1$ integral using
\[ a = \tfrac{d-1}{2}+t+r-2r'-s_2, \qquad b=t+r+r', \qquad c=\tfrac{3t}{2}, \qquad d=s_3, \]
we may evaluate the $s_3$ integral using \eqref{eq:MB1F0Eval}.
The $s_1$ integral becomes a $\pFqName21$ at $\frac{1}{2}$ by \eqref{eq:pFqToMB}, and after a Pfaff transformation \eqref{eq:Pfaff}, we have
\begin{align*}
	\Psi^d=& 2^{4-d-4t-r-r'} \pi^{2-3t} \frac{\Gamma(d-1+2t-r-r') \Gamma(d-1+t+r+r') \Gamma\paren{t+r+r'} \Gamma(\tfrac{t}{2}-r-r')}{\Gamma\paren{d-1+\frac{3t}{2}} \Gamma\paren{\frac{3t}{2}}} \\
	& \times \int_{\mathcal{C}} \Gamma(s_2) \Gamma\paren{d-1+\tfrac{3t}{2}-s_2} \\
	& \times \Biggl(\frac{\Gamma(\frac{d-1}{2}+t+r-2r'-s_2)}{\Gamma(\frac{d+1}{2}+t+r-2r'-s_2)} \pFq21{\frac{d-1}{2}+t+r-2r'-s_2,\tfrac{t}{2}-r-r'}{\frac{d-1}{2}+t+r-2r'-s_2+1}{-1} \\
	& + \frac{\Gamma(\frac{d-1}{2}+t+r-2r'-s_2)}{\Gamma(\frac{d+1}{2}+t+r'-2r-s_2)} \pFq21{\frac{d-1}{2}+t+r'-2r-s_2,\tfrac{t}{2}-r-r'}{\frac{d-1}{2}+t+r'-2r-s_2+1}{-1} \Biggr) \frac{ds_2}{2\pi i}.
\end{align*}

In the second term, we send $s_2 \mapsto d-1+\tfrac{3t}{2}-s_2$, which allows us to use \eqref{eq:2F1toBeta}, so that
\begin{align*}
	\Psi^d=& 2^{4-d-4t-r-r'} \pi^{2-3t} \frac{\Gamma(d-1+2t-r-r') \Gamma(d-1+t+r+r') \Gamma\paren{t+r+r'} \Gamma\paren{\tfrac{t}{2}-r-r'}}{\Gamma\paren{d-1+\frac{3t}{2}} \Gamma\paren{\frac{3t}{2}}} \\
	& \times \int_{\mathcal{C}} \Gamma(s_1)\Gamma\paren{d-1+\tfrac{3t}{2}-s_1} B\paren{\tfrac{d-1}{2}+t+r-2r'-s_1,\tfrac{1-d}{2}-\tfrac{t}{2}+r'-2r+s_1} \frac{ds_1}{2\pi i}.
\end{align*}
The theorem now follows from Barnes' first lemma \eqref{eq:BarnesFirst}.

\section{Kontorovich-Lebedev Inversion}
\label{sect:KontLeb}
If $F(r)$ is holomorphic in a neighborhood $\abs{\Re(r)} < \delta < \frac{1}{10}$, then we can argue that the $Y^+$ integral of $(F^\flat)^\sharp$ converges absolutely (via contour shifting in $r'$ and the Mellin-Barnes integral) and define
\[ F(r,\epsilon) := \int_{Y^+} F^\flat(y) \wbar{\trans{W^{d*}(y,\mu(r))}} (y_1^2 y_2)^\epsilon dy = \int_{\Re(r')=0} F(r') \Psi^d(r',-r,\epsilon) \sinmu^d(r') \, dr', \]
where we assume $\eta := \frac{\delta}{2} > \epsilon > 0$, $\Re(r)=0$ and $r\ne 0$.

Shift the $r'$ integral to $\Re(r')=\eta$, picking up a residue at $r'=\frac{\epsilon}{2}+r$.
The shifted integral is zero in the limit $\epsilon \to 0$ (by the $\Gamma(\frac{3\epsilon}{2})$ in the denominator), and the residue of $\Psi^d$ is
\begin{align*}
	\lim_{\epsilon\to 0^+} \res_{r'=\frac{\epsilon}{2}+r} \Psi^d(r',-r,\epsilon) =& 2^{4-d} \pi^2 \Gamma(d-1) \Gamma(\tfrac{d-1}{2}-3r) \Gamma(\tfrac{d-1}{2}+3r) = \frac{1}{2\pi i \sinmu^d(r)}.
\end{align*}

\section{Asymptotics and Functional Equations of the Whittaker Functions}
As in the previous papers, we will require certain first-term asymptotics of the Whittaker function.
Define $(2d+1)$-dimensional row vectors $\bv^d_j$ with entries $\bv^d_{j,m'} = \delta_{m'=j}$.
Then the asymptotics of the completed Whittaker function are given by the following lemmas.
\begin{lem}
\label{lem:wlAsymps}
	Assume $r \ne 0$, then as $y \to 0$, we have
	\begin{align*}
		W^{d*}(y,r) \sim& \frac{(-1)^d}{\pi} (2\pi)^{\frac{d+3}{2}-3r} p_{\rho+\mu^{w_4}}(y) \Gamma\paren{\tfrac{d-1}{2} + 3r} \bv^d_d \WigDMat{d}(\vpmpm{--} w_l) \\
		&+\frac{1}{\pi} (2\pi)^{\frac{d+3}{2}+3r} p_{\rho+\mu}(y) \Gamma\paren{\tfrac{d-1}{2} - 3r} \bv^d_d \\
		&-(-1)^d \frac{(2\pi)^{\frac{3d+1}{2}+3r} p_{\rho+\mu^{w_3}}(y)}{\Gamma\paren{\frac{d+1}{2}+3r} \sin \pi\paren{\frac{d-1}{2}+3r}} \bv^d_d \WigDMat{d}(\vpmpm{--} w_l) T^d(w_2,\mu^{w_4}),
	\end{align*}
	for $\mu=\mu(r)$, in the sense that $W^{d*}(y,r)$ is a sum of (vector multiples of) three power series with the given leading terms.
\end{lem}

\begin{lem}
\label{lem:w4Asymps}
	Assume $r \ne 0$, then as $y_1 \to 0$, we have
	\begin{align*}
		W^{d*}(y,r) \sim& y_1^{\frac{d+1}{2}-r} y_2^{1-2r} \Lambda^*(r) \Gamma\paren{-\tfrac{d-1}{2}+3r} \frac{(2\pi)^{\frac{3d-1}{2}-3r}}{(d-1)!} \\
	& \qquad \times \paren{i^d \exp\paren{-\tfrac{i\pi}{2}\paren{\tfrac{d-1}{2}-3r}}\bv^d_{-d}+i^{-d} \exp\paren{\tfrac{i\pi}{2}\paren{\tfrac{d-1}{2}-3r}}\bv^d_d} \\
	& \qquad \times T^d(w_3,\mu) \WigDMat{d}(w_2) \mathcal{W}^d\paren{-y_2,\tfrac{d-1}{2}+3r} \\
		&+(-i)^d (2\pi)^{1+6r} y_1^{1+2r} y_2^{\frac{3-d}{2}+r} \Lambda^*(-r) \frac{\Gamma\paren{\frac{d-1}{2}-3r}}{\Gamma\paren{\frac{d+1}{2}+3r}} \bv^d_{-d} \WigDMat{d}(\vpmpm{--}w_l) \\
	& \qquad \times T^d(w_5,\mu^{w_4}) \WigDMat{d}(w_2) \mathcal{W}^d(-y_2,d-1),
	\end{align*}
	for $\mu=\mu(r)$, in the sense that $W^{d*}(y,r)$ is a sum of (vector multiples of) two power series with the given leading terms.
\end{lem}

The functional equation takes the form
\begin{lem}
\label{lem:DualWhitt}
	\begin{align*}
		W^{d*}(g, r) = (-1)^d W^{d*}(\vpmpm{--} g^\iota w_l,-r),
	\end{align*}
	where $g^\iota = w_l \trans{(g^{-1})} w_l$.
\end{lem}
We will use Lemma \ref{lem:wlAsymps} in the proof of Lemma \ref{lem:DualWhitt}, which in turn is used in the proof of Lemma \ref{lem:w4Asymps}.

\subsection{Double asymptotics of the Whittaker Function}
\label{sect:wlWhittAsymps}
We now prove Lemma \ref{lem:wlAsymps}.
Assume $r \ne 0$.
We know that $W^{d*}_{m'}(y,r)$ is a linear combination of power series with leading terms $p_{\rho+\mu^w}(y)$ and it is clear that the terms with $w\in\set{w_2,w_5,w_l}$ do not occur since $W^{d*}_{m'}(y,r) \ll \abs{y_1 y_2}$ for $\Re(r)=0$.
We need to find the coefficients of the remaining first terms; these occur as poles of the integrand in the Mellin-Barnes integral.
As in the definition of the completed Whittaker function, we write $m'=\varepsilon m$ with $\varepsilon=\pm1$ and $m \ge 0$.

The residue at $s_1=-\frac{d-1}{2}+r$, $s_2=2r$ is
\begin{align*}
	R_1 :=& \frac{2}{2^{d+1}\pi} \sqrt{\binom{2d}{d+m}} (2\pi y_1)^{\frac{d+1}{2}-r} (2\pi y_2)^{1-2r} \Gamma\paren{\tfrac{d-1}{2} + 3r}.
\end{align*}
By \cite[(70)-(72)]{HWII}, we have
\begin{align}
\label{eq:WDvmmwl}
	\WigD{d}{d}{m'}(\vpmpm{--} w_l) = (-1)^{m'} \Wigd{d}{d}{m'}(0) = (-1)^{d} 2^{-d} \sqrt{\binom{2d}{d+m'}},
\end{align}
so
\begin{align*}
	R_1 =& \frac{(-1)^d}{\pi} \WigD{d}{d}{m'}(\vpmpm{--} w_l) (2\pi y_1)^{\frac{d+1}{2}-r} (2\pi y_2)^{1-2r} \Gamma\paren{\tfrac{d-1}{2} + 3r}.
\end{align*}

The residue at $s_1=-2r$, $s_2=-\frac{d-1}{2}-r$ is
\begin{align*}
	R_2 :=& \delta_{m=d} \frac{2}{2^{d+1}\pi} (2\pi y_1)^{1+2r} (2\pi y_2)^{\frac{d+1}{2}+r} \Gamma\paren{\tfrac{d-1}{2} - 3r} \sum_{\ell=0}^d \varepsilon^\ell \binom{d}{\ell} \\
	=& \delta_{m'=d} \frac{1}{\pi} (2\pi y_1)^{1+2r} (2\pi y_2)^{\frac{d+1}{2}+r} \Gamma\paren{\tfrac{d-1}{2} - 3r}.
\end{align*}

The residue at $s_1=-\frac{d-1}{2}+r$, $s_2=-\frac{d-1}{2}-r$ is
\begin{align*}
	R_3 :=& \frac{1}{2^{d+1}\pi} \sqrt{\binom{2d}{d+m}} (2\pi y_1)^{\frac{d+1}{2}-r} (2\pi y_2)^{\frac{d+1}{2}+r} \Gamma\paren{\tfrac{\frac{d+1}{2}-m+3r}{2}} \sum_{\ell=0}^{m} \varepsilon^\ell \binom{m}{\ell} \frac{\Gamma\paren{\frac{\ell-\frac{d-1}{2}-3r}{2}}}{\Gamma\paren{\frac{1-m+\ell}{2}}}.
\end{align*}
The terms with $\ell \not\equiv m\pmod{2}$ are zero, and the sum of the remaining terms may be evaluated by converting to a $\pFqName21$ at 1 and applying \eqref{eq:GaussThm}:
\begin{align*}
	\sum_{\ell=0}^{m} \varepsilon^\ell \binom{m}{\ell} \frac{\Gamma\paren{\frac{\ell-\frac{d-1}{2}-3r}{2}}}{\Gamma\paren{\frac{1-m+\ell}{2}}} =& -\varepsilon^m \frac{2^{\frac{d-1}{2}+3r}}{\Gamma\paren{\frac{d+1}{2}+3r}} \frac{\Gamma\paren{\frac{\frac{d+1}{2}+3r+m}{2}}}{\sin\frac{\pi}{2}\paren{\frac{d-1}{2}+3r-m}}.
\end{align*}
Using \eqref{eq:Tdw2}-\eqref{eq:WhittGammas}, we may write this as
\begin{align*}
	R_3 =& -(-1)^d (2\pi)^{\frac{d-1}{2}+3r} \frac{(2\pi y_1)^{\frac{d+1}{2}-r} (2\pi y_2)^{\frac{d+1}{2}+r}}{\Gamma\paren{\frac{d+1}{2}+3r} \sin \pi\paren{\frac{d-1}{2}+3r}} \paren{\WigDMat{d}(\vpmpm{--} w_l) T^d(w_2,\mu(r)^{w_4})}_{d,m'}.
\end{align*}

Then in a formal sense, we have
\[ W^{d*}_{m'}(y,r) \sim R_1+R_2+R_3 \]
as $y \to 0$, or in other words, $W^{d*}_{m'}$ is given by a sum of three power series with those leading terms.
This will be sufficient to identify the particular linear combination of power series occuring in the $K^d_{w_l}$ functions, as those power series are also distinguished by their leading terms.

\subsection{The dual Whittaker function}
\label{sect:dualWhitt}
We must briefly resort to the differential operators $Y^a$.
As in \cite[sect. 6.3]{HWII}, if we take the dual $\wcheck{f}(g) := f(\vpmpm{--} g^\iota w_l)$ of the function $f(g) :=W^{d*}(g, -r)$, then the action of the lowering operators is
\[ Y^{-1} \wcheck{f}(g) = -\wcheck{Y^{-1}f}(g) = 0, \qquad Y^{-2} \wcheck{f}(g) = -\wcheck{Y^{-2}f}(g) = 0, \]
by the minimality of $f$.

Assume $\Re(r)=0$, then from \cite[sect. 6.3]{HWII} and \eqref{eq:WhittFEs}, we have
\[ \wcheck{f}(g) = \paren{\Lambda^*(-r) \bv^d_{-d} \WigDMat{d}(\vpmpm{--} w_l) T^d(w_2,\mu(r)^{w_4})} W^d(g,\mu(r)^{w_3},\psi_{1,1}), \]
so $\wcheck{f}(g)$ lies in the rowspace of $W^d(g,\mu(r)^{w_3},\psi_{1,1})$.
If we also assume $r \ne 0$, then \cite[prop. 17]{HWII}, \cite[prop. 15]{HWII}, and \eqref{eq:WhittCompletion} (and the $w_3$ functional equation \eqref{eq:WhittFEs}) imply $\wcheck{f}(g) = C W^{d*}(g, r)$ for some scalar $C=C(d,r)$.
But then Lemma \ref{lem:wlAsymps} implies $C=(-1)^d$ by comparing asymptotics, and this extends to an equality of meromorphic functions.

\subsection{Single asymptotics of the Whittaker Function}
\label{sect:w4WhittAsymps}
We now prove Lemma \ref{lem:w4Asymps}.
If $\Re(r)$ is large, then as in \cite[sect. 7.2]{HWII}, we have
\begin{align*}
	W^{d*}(y,r) \sim& p_{\rho+\mu^{w_l}}(y) \Lambda^*(r) \frac{(-2\pi)^d y_1^{d-1}}{(d-1)!} \bv^d_{-d} \int_{\R^2} (1+u_3^2)^{\frac{-1+\frac{d-1}{2}-3r}{2}} (1+u_2^2)^{\frac{-1-\frac{d-1}{2}-3r}{2}} \\
	& \qquad \times \WigDMat{d}(w_3) \Dtildek{d}{\frac{-u_3-i}{\sqrt{1+u_3^2}}} \WigDMat{d}(w_3) \Dtildek{d}{\frac{1-iu_2}{\sqrt{1+u_2^2}}} \e{-y_2 u_2} du,
\end{align*}
as $y_1 \to 0$.
The integrals may then be expressed in terms of the $\mathcal{W}^d$ function as in \eqref{eq:classWhittDef}:
\begin{align*}
	W^{d*}(y,r) \sim& y_1^{\frac{d+1}{2}-r} y_2^{1-2r} \Lambda^*(r) \frac{(2\pi)^d}{(d-1)!} \bv^d_{-d} \WigDMat{d}(w_3) \mathcal{W}^d\paren{0,-\tfrac{d-1}{2}+3r} \\
	& \qquad \times \WigDMat{d}(w_5) \mathcal{W}^d\paren{-y_2,\tfrac{d-1}{2}+3r}.
\end{align*}

From the functional equation of Lemma \ref{lem:DualWhitt}, for $\Re(r)$ highly negative, we have
\begin{align*}
	W^{d*}(y,r) \sim& (-1)^d y_1^{1+2r} y_2^{\frac{3-d}{2}+r} \Lambda^*(-r) \bv^d_{-d} W^d(I,-\mu^{w_2},\psi_{y_2,0}) \WigDMat{d}(\vpmpm{--}w_l),
\end{align*}
as $y_1 \to 0$.
We then insert the computed value \cite[(3.11)]{HWI}:
\begin{align*}
	W^{d*}(y,r) \sim& (-1)^d y_1^{1+2r} y_2^{\frac{3-d}{2}+r} \Lambda^*(-r) \bv^d_{-d} \WigDMat{d}(\vpmpm{--}w_l) \mathcal{W}^d\paren{0,-\tfrac{d-1}{2}-3r} \\
	& \qquad \times \WigDMat{d}(w_3) \mathcal{W}^d\paren{0,\tfrac{d-1}{2}-3r} \WigDMat{d}(w_5) \mathcal{W}^d(-y_2,d-1).
\end{align*}

As with the $y \to 0$ asymptotic, we know that the term $y_1^{\frac{3-d}{2}-r}$ does not occur.
Applying \eqref{eq:W0toGammaW}, \eqref{eq:WhittFEvsV} and
\begin{align}
\label{eq:bvddvw2}
	\bv^d_{\pm d} \WigDMat{d}(v_{\varepsilon_1,\varepsilon_2}) = (\varepsilon_1 \varepsilon_2)^d \bv^d_{\pm\varepsilon_2 d}, \qquad \bv^d_{\pm d} \WigDMat{d}(w_2) = (-1)^d i^{\pm d} \bv^d_{\mp d}
\end{align}
completes the lemma, keeping in mind that $\WigDMat{d}(\vpmpm{+-}w_2)=\Dtildek{d}{i}$ commutes with diagonal matrices such as $T^d(w_2,\mu)$.

\section{Kuznetsov's Formula}
\label{sect:KuzForm}
We consider a Fourier coefficient of a Poincar\'e series of the form
\begin{align}
\label{eq:MainPoincare}
	P_m(g,F) =& \sum_{\gamma\in U(\Z)\backslash\Gamma} \int_{\Re(r)=0} F(r) W^{d*}(\wtilde{m} \gamma g,r) \sinmu^d(r)\,dr, \qquad \wtilde{m}=\diag(m_1 m_2,m_1,1).
\end{align}
There are some technical issues with the convergence of this series, especially for $d=3$ and $d=2$, and we will discuss them in section \ref{sect:Weight2}.

Define the Fourier-Whittaker coefficients of a Maass form $\xi$ with Langlands parameters $\mu(r_\xi)$ by
\begin{align}
\label{eq:FWcoefDef}
	\int_{U(\Z)\backslash U(\R)} \xi(xyk) \wbar{\psi_m(x)} dx = \frac{\rho_\xi^*(n)}{\abs{m_1 m_2}} W^{d*}(\wtilde{m} yk, r_\xi),
\end{align}
and define the integral transform
\begin{align}
\label{eq:HwTildeDef}
	\wtilde{H}_w(F; y, g) =& \frac{1}{\abs{y_1 y_2}} \int_{\wbar{U}_w(\R)} \int_{\Re(r)=0} F(r) W^{d*}(ywxg,r) \sinmu^d(r)\,dr \,\wbar{\psi_{1,1}(x)} dx,
\end{align}
for $w\in W$, $y\in Y$, $g\in G$.

The spectral expansion and Bruhat decomposition give the pre-Kuznetsov formula
\begin{align}
\label{eq:PreKuzSpectral}
	& \int_\mathcal{B} F(\mu_\xi) \wbar{\rho_\xi^*(m)} \rho_\xi^*(n) W^{d*}(\wtilde{n} yk, \mu_\xi) \, d\xi \\
	&= \abs{\frac{n_1 n_2}{m_1 m_2}} \int_{U(\Z)\backslash U(\R)} P_m(xyk,F) \wbar{\psi_n(x)} dx \nonumber \\
	&= \sum_{w\in W} \sum_{v\in V} \sum_{c_1,c_2\ge1} \frac{S_w(\psi_m,\psi_n^v,c)}{c_1 c_2} \wtilde{H}_w\paren{F; \wtilde{m} c w v \wtilde{n}^{-1} w^{-1},\wtilde{n} yk}, \nonumber
\end{align}
where $\int_\mathcal{B} \ldots d\xi$ serves as an abbreviation for the sums and integrals occuring in the spectral expansion \eqref{eq:SpectralExpand}.
The details of the Bruhat decomposition can be found in \cite[sects 2.2.3,2.2.5]{MeThesis}.

Note: To be precise, in the development of the Kuznetsov formula, we must initially require $F$ to be holomorphic on $\abs{\Re(r)} < \frac{d}{6}+\delta$ and have sufficient exponential decay to overcome the growth of the Fourier-Whittaker coefficients for absolute convergence of the Poincar\'e series.
We may relax to $\abs{\Re(r)} < \frac{1}{4}+\delta$ once we reach the pre-Kuznetsov formula.

In section \ref{sect:KuzKernelsProof} below, we show
\begin{lem}
\label{lem:KuzKernels}
	Let $F$ be holomorphic and Schwartz-class on a neigborhood of $\Re(\mu)=0$, then for $w=I,w_4,w_5,w_l$, we have
	\begin{align*}
		\wtilde{H}_w(F; y, g) =& \frac{1}{\abs{y_1 y_2}} \int_{\Re(r)=0} F(r) K_w^d(y,r) W^{d*}(g,r) \sinmu^d(r)\,dr,
	\end{align*}
	with $K_w^d(y,r)$ as in \eqref{eq:KwIEval}-\eqref{eq:KwlEval}.
\end{lem}

Then replacing $F(r)$ with
\[ F(r) \trans{\wbar{W^{d*}(\wtilde{n} yk,r)}} (y_1^2 y_2) \]
(in a suitable manner), and integrating in $y$ with Stade's formula gives the theorem.
There is a small technical point that the pre-Kuznetsov formula \eqref{eq:PreKuzSpectral} is an equality of vectors and we wish to apply Stade's formula, which involves a dot product, \textit{inside} the $r$-integral; this may be accomplished, e.g., by taking the central entry of \eqref{eq:PreKuzSpectral}, replacing $F(r)$ with the central entry of the previous display, and integrating over $k$.

\subsection{Power series for the Kuznetsov kernel functions}
The functions $K_w^d(y,r)$ are solutions to the differential equations
\[ \Delta_i K_w(g,r) = \lambda_i(\mu(r)) K_w(g,r), \qquad K_w(u g (wu'w^{-1}),r) = \psi_{1,1}(uu') K_w(g,r), \]
where $g\in G$, $w\in W$, $u\in U(\R)$, $u'\in U_w(\R)$ and
\[ \lambda_1(\mu) = 1-\tfrac{\mu_1^2+\mu_2^2+\mu_3^2}{2}, \qquad \lambda_2(\mu) = \mu_1 \mu_2 \mu_3. \]
These were solved in the paper \cite{SpectralKuz}, under the assumption that $\mu_i-\mu_j \notin \Z$, $i\ne j$, but of course that fails for $\mu=\mu(r)$.

When $w=w_l$, the original power-series solutions are
\begin{align}
\label{eq:JwlDef}
	J_{w_l}(y,\mu) =& \abs{4\pi^2 y_1}^{1-\mu_3} \abs{4\pi^2 y_2}^{1+\mu_1} \sum_{n_1,n_2\ge 0} \frac{\Gamma\paren{n_1+n_2+\mu_1-\mu_3+1} \, (4\pi^2 y_1)^{n_1} (4\pi^2 y_2)^{n_2}}{\prod_{i=1}^3 \Gamma\paren{n_1+\mu_i-\mu_3+1}\Gamma\paren{n_2+\mu_1-\mu_i+1}}.
\end{align}
By comparing asymptotics, we can see that $J_{w_l}(y,\mu)$, $J_{w_l}(y,\mu^{w_3})$ and $J_{w_l}(y,\mu^{w_4})$ are distinct solutions, but we have the relations
\begin{align*}
	J_{w_l}(y,\mu(r)^{w_2}) =& \sgn(y_2)^{d-1} J_{w_l}(y,\mu(r)), \\
	J_{w_l}(y,\mu(r)^{w_l}) =& \sgn(y_1)^{d-1} J_{w_l}(y,\mu(r)^{w_4}), \\
	J_{w_l}(y,\mu(r)^{w_5}) =& \sgn(y_1 y_2)^{d-1} J_{w_l}(y,\mu(r)^{w_3}).
\end{align*}
For $w=w_4,w_l$ and $w',\tilde{w}\in W$, we define the linear combination
\[ Y^d_w(y,\mu,w',\tilde{w},\alpha) = \frac{J_w(y,\mu^{\tilde{w}})-\alpha^{d-1} J_w(y,\mu^{w'})}{\sin\pi(\mu_1-\mu_2)}. \]
Then the remaining three long-element solutions and their first-term asymptotics as $y \to 0$ are
\begin{align*}
	Y^d_{w_l}(y,\mu(r), I, w_2,\sgn(y_2)) \sim& \frac{\abs{4\pi^2 y_1}^{1+2r}\abs{4\pi^2 y_2}^{\frac{3-d}{2}+r}(d-2)!}{\pi \Gamma\paren{\frac{3-d}{2}+3r} \Gamma\paren{\frac{d+1}{2}+3r}}, \\
	Y^d_{w_l}(y,\mu(r),w_3, w_5,\sgn(y_1 y_2)) \sim& \frac{\abs{4\pi^2 y_1}^{\frac{3-d}{2}-r}\abs{4\pi^2 y_2}^{\frac{3-d}{2}+r}(d-2)!}{\pi \Gamma\paren{\frac{3-d}{2}+3r} \Gamma\paren{\frac{3-d}{2}-3r}}, \\
	Y^d_{w_l}(y,\mu(r),w_4, w_l,\sgn(y_1)) \sim& \frac{\abs{4\pi^2 y_1}^{\frac{3-d}{2}-r}\abs{4\pi^2 y_2}^{1-2r}(d-2)!}{\pi \Gamma\paren{\frac{3-d}{2}-3r} \Gamma\paren{\frac{d+1}{2}-3r}}.
\end{align*}
In each case, either $y_1$ or $y_2$ has an exponent less than one, and will not appear in the Kuznetsov kernel functions.
(For $d \le 3$, the main asymptotics as $y \to 0$ are actually given by some logarithmic factors; again, such terms cannot appear in the kernel functions.)

For $w=w_4$, the power-series solutions are (see \cite[sect. 4.5]{WeylI}),
\begin{align}
\label{eq:Jw4Def}
	J_{w_4}(y,\mu) =& \abs{8\pi^3 y_1}^{1-\mu_3} \sum_{n=0}^\infty \frac{(-8\pi^3 i y_1)^n}{n! \, \Gamma\paren{n+1+\mu_1-\mu_3} \, \Gamma\paren{n+1+\mu_2-\mu_3}},
\end{align}
and $J_{w_4}(y,\mu(r)^{w_4})$ is distinct from $J_{w_4}(y,\mu(r))$, but
\begin{align*}
	J_{w_4}(y,\mu(r)^{w_5}) =& (i \sgn(y_1))^{d-1} J_{w_4}(y,\mu(r)^{w_4}).
\end{align*}
The requisite third solution is given by
\begin{align*}
	Y^d_{w_4}(y,\mu(r), w_4,w_5,i\sgn(y_1)) \sim& \frac{\abs{8\pi^3 y_1}^{\frac{3-d}{2}-r}(d-2)!}{\pi \Gamma\paren{\frac{3-d}{2}-3r}},
\end{align*}
and again, this will not appear in the Kuznetsov kernel functions.

\subsection{The weight functions}
\label{sect:KuzKernelsProof}
Having done the relevant technical work related to the analytic continuation in previous papers, we regard the functions $K_w^d(y,r)$ as being defined by the Riemann integral
\begin{align*}
	K_w^d(y,r) W^{d*}(g,r) =& \int_{\wbar{U}_w(\R)} W^{d*}(ywxg,r) \,\wbar{\psi_{1,1}(x)} dx,
\end{align*}

Let $x^* y^* k^* = wxg$ and replace $y \mapsto v_{\varepsilon_1,\varepsilon_2} y$ with $y\in Y^+$, then formally
\begin{align*}
	K_w^d(y,r) W^{d*}(g,r) =& \int_{\wbar{U}_w(\R)} W^{d*}(y y^*,r) \WigDMat{d}(v_{\varepsilon_1,\varepsilon_2} k^*) \,\psi_{v_{\varepsilon_1,\varepsilon_2} y}(x^*) \wbar{\psi_{1,1}(x)} dx.
\end{align*}

The process for obtaining $K_w^d(y,r)$ is the same as for the previous two cases, but we must shift the $s$-integrals to $\Re(s)=(-\frac{d-1}{2}-\epsilon,-\frac{d-1}{2}-\epsilon)$ to see the term $y_1^{\frac{d+1}{2}-r} y_2^{\frac{d+1}{2}+r}$, and the test function needs holomorphy out to $\abs{\Re(r)} < \frac{d}{6}+\epsilon$ for absolute convergence of the integrals and sums of Kloosterman sums for the terms $y_1^{1+2r} y_2^{\frac{d+1}{2}+r}$ and $y_1^{\frac{d+1}{2}-r} y_2^{1-2r}$.
Lastly, the latter two terms require many more applications of the integration by parts proceedure on the function $X_3'$ in \cite[sect. 2.6.2]{SpectralKuz} to reach $\Re(r)=0$, but as explained in the preceeding paper, this is always possible.

\subsubsection{The long element function}
As $y \to 0$,
\begin{align*}
	& K_{w_l}^d(v_{\varepsilon_1,\varepsilon_2} y,r) W^{d*}(g,r) \\
	&\sim 2 (-1)^d p_{\rho+\mu^{w_4}}(y) (2\pi)^{\frac{d+1}{2}-3r} \Gamma\paren{\tfrac{d-1}{2} + 3r} \bv^d_d \WigDMat{d}(\vpmpm{--} w_l v_{\varepsilon_1,\varepsilon_2}) \\
	& \qquad \times \int_{U(\R)} p_{\rho+\mu^{w_4}}(y^*) \WigDMat{d}(k^*) \wbar{\psi_{1,1}(x)} dx \\
	& + 2 p_{\rho+\mu}(y) (2\pi)^{\frac{d+1}{2}+3r} \Gamma\paren{\tfrac{d-1}{2} - 3r} \bv^d_d \WigDMat{d}(v_{\varepsilon_1,\varepsilon_2}) \int_{U(\R)} p_{\rho+\mu}(y^*) \WigDMat{d}(k^*) \wbar{\psi_{1,1}(x)} dx \\
	&-(-1)^d \frac{(2\pi)^{\frac{3d+1}{2}+3r} p_{\rho+\mu^{w_3}}(y)}{\Gamma\paren{\frac{d+1}{2}+3r} \sin \pi\paren{\frac{d-1}{2}+3r}} \bv^d_d \WigDMat{d}(\vpmpm{--} w_l) T^d(w_2,\mu^{w_4}) \WigDMat{d}(v_{\varepsilon_1,\varepsilon_2}) \\
	& \qquad \times \int_{U(\R)} p_{\rho+\mu^{w_3}}(y^*) \WigDMat{d}(k^*) \wbar{\psi_{1,1}(x)} dx,
\end{align*}
in the sense of Lemma \ref{lem:wlAsymps}.

The $x$ integrals give incomplete Whittaker functions by comparing to the definition \eqref{eq:JacWhittDef}, and the functional equations \eqref{eq:WhittFEs} and Lemma \ref{lem:DualWhitt} (keeping in mind \eqref{eq:WhittFEvsV} and \eqref{eq:bvddvw2}) imply
\begin{align*}
	K_{w_l}^d(v_{\varepsilon_1,\varepsilon_2} y,r) \sim& 2 (-\varepsilon_2)^d p_{\rho+\mu^{w_4}}(y) \delta_{\varepsilon_1=-1} \frac{(2\pi)^{\frac{d+1}{2}-3r} \Gamma\paren{\tfrac{d-1}{2} + 3r}}{\Lambda^*(-r)} \\
	& + 2 (-\varepsilon_1)^d p_{\rho+\mu}(y) \delta_{\varepsilon_2=-1} \frac{(2\pi)^{\frac{d+1}{2}+3r} \Gamma\paren{\tfrac{d-1}{2} - 3r}}{\Lambda^*(r)} \\
	&-(\varepsilon_1)^d p_{\rho+\mu^{w_3}}(y) \delta_{\varepsilon_1 \varepsilon_2=-1} \frac{(2\pi)^{\frac{3d+1}{2}+3r}}{\Gamma\paren{\frac{d+1}{2}+3r} \sin \pi\paren{\frac{d-1}{2}+3r} \Lambda^*(-r)}.
\end{align*}
Then comparing asymptotics with $J_{w_l}(y,\mu(r)^{w_4})$, $J_{w_l}(y,\mu(r))$, and $J_{w_l}(y,\mu(r)^{w_3})$ gives \eqref{eq:KwlEval}.

\subsubsection{The $w_4$ function}
As $y_1 \to 0$ on $\varepsilon_2=y_2=1$, with $\mu=\mu(r)$,
\begin{align*}
	& K_{w_4}^d(\vpmpm{\varepsilon_1,1}y,r) W^{d*}(g,r) \\
	&\sim y_1^{\frac{d+1}{2}-r} \Lambda^*(r) \Gamma\paren{-\tfrac{d-1}{2}+3r} \frac{(2\pi)^{\frac{3d-1}{2}-3r}}{(d-1)!} \\
& \qquad \times \paren{i^d \exp\paren{-\tfrac{i\pi}{2}\paren{\tfrac{d-1}{2}-3r}}\bv^d_{-d}+i^{-d} \exp\paren{\tfrac{i\pi}{2}\paren{\tfrac{d-1}{2}-3r}}\bv^d_d} T^d(w_3,\mu) \WigDMat{d}(w_2 v_{\varepsilon_1,1}) \\
	& \qquad \times \int_{\wbar{U}_{w_4}(\R)} (y_1^*)^{1-\mu_2} (y_2^*)^{1+\mu_3} \mathcal{W}^d\paren{-y_2^*,\mu_1-\mu_3} \WigDMat{d}(k^*) \,\e{x_2^*-x_2} dx \\
	&\quad +(-i)^d (2\pi)^{1+6r} y_1^{1+2r} \Lambda^*(-r) \frac{\Gamma\paren{\frac{d-1}{2}-3r}}{\Gamma\paren{\frac{d+1}{2}+3r}} \bv^d_{-d} \WigDMat{d}(\vpmpm{--}w_l) T^d(w_5,\mu^{w_4}) \WigDMat{d}(w_2 v_{\varepsilon_1,1}) \\
	& \qquad \times \int_{\wbar{U}_{w_4}(\R)} (y_1^*)^{1-\mu_3} (y_2^*)^{1+\mu_2} \mathcal{W}^d(-y_2^*,\mu_1-\mu_2) \WigDMat{d}(k^*) \,\e{x_2^*-x_2} dx,
\end{align*}
in the sense of Lemma \ref{lem:w4Asymps}.

The $x$ integrals give incomplete Whittaker functions by comparing to \cite[(3.22)]{HWI}, so that
\begin{align*}
	& K_{w_4}^d(\vpmpm{\varepsilon_1,1}y,r) W^{d*}(g,r) \\
	&\sim y_1^{\frac{d+1}{2}-r} \Lambda^*(r) \Gamma\paren{-\tfrac{d-1}{2}+3r} \frac{(2\pi)^{\frac{3d-1}{2}-3r}}{(d-1)!} \\
& \qquad \times \paren{i^d \exp\paren{-\tfrac{i\pi}{2}\paren{\tfrac{d-1}{2}-3r}}\bv^d_{-d}+i^{-d} \exp\paren{\tfrac{i\pi}{2}\paren{\tfrac{d-1}{2}-3r}}\bv^d_d} \\
& \qquad \times \WigDMat{d}(v_{\varepsilon_1,\varepsilon_1}) T^d(w_3,\mu) W^d(g,\mu^{w_3},\psi_{1,1}) \\
	&\quad +(-i)^d (2\pi)^{1+6r} y_1^{1+2r} \Lambda^*(-r) \frac{\Gamma\paren{\frac{d-1}{2}-3r}}{\Gamma\paren{\frac{d+1}{2}+3r}} \bv^d_{-d} \WigDMat{d}(v_{\varepsilon_1,1}) \WigDMat{d}(\vpmpm{--}w_l) \\
	& \qquad \times T^d(w_5,\mu^{w_4}) W^d(g,\mu,\psi_{1,1}).
\end{align*}

Then applying the functional equations \eqref{eq:WhittFEs} and Lemma \ref{lem:DualWhitt}, this becomes
\begin{align*}
	K_{w_4}^d(\vpmpm{\varepsilon_1,1}y,r) \sim& y_1^{\frac{d+1}{2}-r} \Gamma\paren{-\tfrac{d-1}{2}+3r} \frac{(2\pi)^{\frac{3d-1}{2}-3r}}{(d-1)!} (\varepsilon_1 i)^d \exp\paren{-\varepsilon_1 \tfrac{i\pi}{2}\paren{\tfrac{d-1}{2}-3r}} \\
	&+(\varepsilon_1 i)^d (2\pi)^{1+6r} y_1^{1+2r} \frac{\Gamma\paren{\frac{d-1}{2}-3r}}{\Gamma\paren{\frac{d+1}{2}+3r}}.
\end{align*}

The expression \eqref{eq:Kw4Eval} follows by comparing asymptotics with $J_{w_4}(y,\mu(r))$ and $J_{w_4}(y,\mu(r)^{w_4})$.

\subsubsection{The $w_5$ function}
The definition of the involution $\iota$ and two applications of Lemma \ref{lem:DualWhitt} give
\begin{align*}
	K_{w_5}^d(y,r) W^{d*}(g,r) =& (-1)^d \int_{\wbar{U}_{w_4}(\R)} W^{d*}(\vpmpm{--} y^\iota w_4 \vpmpm{--}x \vpmpm{--}g^\iota w_l,-r) \,\wbar{\psi_{1,1}(x)} dx \\
	=& K_{w_4}^d(y^\iota \vpmpm{-+};-r) W^{d*}(g,r),
\end{align*}
and this implies \eqref{eq:Kw5}.

\section{The Technical Weyl Law}
\label{sect:TechnicalWeylLaw}
We now prove Theorem \ref{thm:TechnicalWeyl}.
By either Stirling's formula or the Phragm\'en-Lindel\"of principle, when $\Re(s) > d$ is in some fixed compact set, we have
\begin{align*}
	\abs{Q(d,s)} \ll& \paren{d+\abs{s}}^{2\Re(s)-1}.
\end{align*}
(Note that, eg. $\abs{Q(d,s)}=\abs{\frac{d-1}{2}+s}^{-1}$ on $\Re(s)=0$.)
Stirling's formula also implies
\begin{align*}
	\abs{B^{\varepsilon}_{w_l}(s,r)} \ll& \paren{1+\abs{s_1+3r}}^{-1-\epsilon} \paren{1+\abs{s_2-3r}}^{-1-\epsilon} \paren{1+\abs{s_1+s_2}}^{\frac{3}{2}+2\epsilon},
\end{align*}
on $\Re(r)=0$, $\Re(s)=\paren{-\frac{1}{2}-\epsilon,-\frac{1}{2}-\epsilon}$.

On $\Re(r)=0$, as in \cite[lem. 15]{WeylI}, we can show
\begin{align*}
	& \int_{\Re(s)=\paren{-\frac{1}{2}-\epsilon,-\frac{1}{2}-\epsilon}} \frac{\paren{1+\abs{s_1+s_2}}^{\frac{3}{2}+2\epsilon}}{\paren{1+\abs{s_1+r}}^{1+\epsilon} \paren{1+\abs{s_2-r}}^{1+\epsilon}} \paren{d+\abs{s_1}}^{-2} \paren{d+\abs{s_2}}^{-2} \abs{ds_1 ds_2} \\
	&\ll \frac{\Min{d,\abs{r}}^{\frac{1}{2}+\epsilon}}{d(d+\abs{r})^2}
\end{align*}
by following the methods of \cite[lem. 4 and 6]{Me01}.

Starting from \eqref{eq:HwDef} and the Mellin-Barnes integral \eqref{eq:KwlMB}, we shift the $s$-contours back to $\Re(s)=\paren{-\frac{1}{2}-\eta,-\frac{1}{2}-\eta}$, picking up poles at $s_1=-3r$ or $s_2=3r$ (but not both simultaneously), where we shift the $r$ contour to $\pm(\frac{1}{4}+\eta)$ and place the remaining $s$-contour at $-\frac{1}{4}-\eta$, giving
\begin{align*}
	& \abs{y_1 y_2}^{-\frac{1}{2}-\epsilon} \abs{H_{w_l}(F;y)} \\
	&\ll \int_{\Re(r)=0} \abs{F(r)} \int_{\Re(s)=\paren{-\frac{1}{2}-\eta,-\frac{1}{2}-\eta}} \paren{1+\abs{s_1+3r}}^{-1-\epsilon} \paren{1+\abs{s_2-3r}}^{-1-\epsilon} \\
	&\qquad \times \paren{1+\abs{s_1+s_2}}^{\frac{3}{2}+2\epsilon} \paren{d+\abs{s_1}}^{-2} \paren{d+\abs{s_2}}^{-2} \abs{ds_1 ds_2} d(d+\abs{r})^2 \abs{dr} \\
	&+ \int_{\Re(r)=-\frac{1}{4}-\eta} \abs{F(r)} \int_{\Re(s_1)=-\frac{1}{4}-\eta} \paren{d+\abs{s_1}}^{-\frac{3}{2}} \abs{ds_1} \abs{Q(d,3r)} d(d+\abs{r})^2 \abs{dr} \\
	&+ \int_{\Re(r)=\frac{1}{4}+\eta} \abs{F(r)} \int_{\Re(s_2)=-\frac{1}{4}-\eta} \paren{d+\abs{s_2}}^{-\frac{3}{2}} \abs{ds_2} \abs{Q(d,-3r)} d(d+\abs{r})^2 \abs{dr} \\
	&\ll E_1 +E_2.
\end{align*}
Note: The residues in $s$ are actually given by $J$-Bessel functions, so we could apply known bounds for those, but the bounds above are sufficient for our purposes.

For the $w_4$ term, we start with the second form of \eqref{eq:Kw4MB}, and shift to $\Re(s) = -\frac{1}{2}$.
(We only need $\Re(s)=-\epsilon$, but shifting farther would give a better bound; this choice gives the more concise statement).
The residue at $s=-3r$ we shift up to $\Re(r)=\frac{1}{4}$.
The residue is trivial to handle, and for the shifted contour, we use
\begin{align*}
	\int_{\Re(s)=0} \paren{1+\abs{s+r}}^{-1} \paren{d+\abs{s}}^{-2} \abs{ds} \ll& \frac{1}{d(d+\abs{r})}.
\end{align*}
The $w_5$ term is handled by symmetry.

For the Eisenstein series term, we assume $\mathcal{S}^{d*}_2$ contains Hecke eigenforms so that we may (skip ahead a little and) use $\mathcal{E}$ in the form from Theorem \ref{thm:HeckeKuznetsov}.
It is well-known that $\abs{\mathcal{S}^{d*}_2} \ll d$ and the quotient by $L$-functions is bounded by $d^\epsilon (1+\abs{r})^\epsilon$ \cite{HoffLock,HoffRam} (see the second remark on page 164 of \cite{HoffLock}).

\section{The Weyl Law}
\label{sect:WeylLaw}
In this section, we prove Corollary \ref{cor:Weyl}.
Taking a test function $F(r) = \exp (r-Tr')^2$, it follows from Theorem \ref{thm:TechnicalWeyl} that
\begin{align}
\label{eq:BoundaryBound}
	\sum_{\abs{r_\varphi-T r'} < 100} \frac{\abs{\rho_\varphi^*(1)}^2}{\cosmu^d(r_\varphi)} \ll d(d+T)^2
\end{align}
for $r'\in i\R$.

Let $\chi_{\abs{r_\varphi-T r'} < M}$ be the characteristic function of the set $\abs{r_\varphi-T r'} < M$, then we define our test function by convolution with an approximation to the identity:
\[ F(r) = -i\sqrt{\frac{\log (d+T)}{\pi}} \int_{\Re(r')=0} \chi_{\abs{r_\varphi-T r'} < M}(r-r') (d+T)^{(r'^2)} dr', \]
Substituting $r'\mapsto r-r'$, this extends to an entire function of $r$.

As in the previous paper,
\begin{align}
\label{eq:TestFunOnSymmBound}
\begin{aligned}
	0 < F(r) < 1 && \text{ on } & r\in i\R \\
	\abs{\chi_{\abs{r_\varphi-T r'} < M}(r)-F(r)} \ll (d+T)^{-100} && \text{ on } & r\in i\R, \abs{r-T r'\pm iM} \ge 10,
\end{aligned}
\end{align}
and in general $F$ satisfies the bound
\begin{align}
\label{eq:TestFunOffSymmBound}
	F(r) \ll (d+T)^{\Re(r)^2+\epsilon} \chi_{\abs{r_\varphi-T r'} < M+10}(i\Im(r))+(\abs{r}+d+T)^{-97} && \text{ on } & \abs{\Re(r)}<1.
\end{align}
Applying \eqref{eq:TestFunOffSymmBound} in Theorem \ref{thm:TechnicalWeyl}, we see that the integrals $E_1$ and $E_2$ are small compared to the error $d(d+T)^2$ resulting from the sharp cut-off.
That error is again obtained by covering the inflated boundary $\abs{r-Tr'\pm iM} < 10$ by 2 balls of radius 11, and applying \eqref{eq:BoundaryBound} and \eqref{eq:TestFunOnSymmBound}.

\section{Rankin-Selberg}
\label{sect:RS}
The computation on $\varphi\in\mathcal{S}^{d*}_3$ in Theorem \ref{thm:HeckeKuznetsov} follows precisely as in \cite[sect. 9.2]{WeylI}, but some more work is required for the maximal parabolic Eisenstein series:
From \cite[(3.31),(5.19)]{HWI}, \eqref{eq:WhittCompletion}, and
\begin{align*}
	\Lambda^*(r)\wbar{\Lambda^*(r)} =& \frac{(d-1)!}{\pi^d \cosmu^d(r)}
\end{align*}
we have
\begin{align}
	\frac{\wbar{\rho_\phi^*(m;r)} \rho_\phi^*(n;r)}{\cosmu^d(r)} =& 4\frac{\pi^d}{(d-1)!} \frac{\wbar{\lambda_\phi(m,r)} \lambda_\phi(n,r)}{L(\phi,1+3r) L(\phi,1-3r)},
\end{align}
when $\phi$ is Hecke-normalized (as in the $\Phi^d_{\what{H}}$ normalization of \cite[sect. 5.3]{HWI}).
When $\phi$ is $L^2$-normalized, this becomes
\begin{align}
\label{eq:MaxParaRS}
	\frac{\wbar{\rho_\phi^*(m;r)} \rho_\phi^*(n;r)}{\cosmu^d(r)} =& 2\pi \frac{\wbar{\lambda_\phi(m,r)} \lambda_\phi(n,r)}{L(\phi,1+3r) L(\phi,1-3r) L(1,\AdSq \phi)}.
\end{align}

\section{Absolute convergence and weight $d=2$}
\label{sect:Weight2}
In this section, we give an upper bound for the Weyl law at $d=2$, as the naive Kuznetsov formula just fails to converge absolutely.
Before we begin, we take a moment to discuss the convergence of the Kuznetsov formula for all $d \ge 3$.

In our development of the spectral Kuznetsov formulae above, we are using implicitly the polynomial dependence on $T$ in the Weyl law.
Ideally, this would follow from M\"uller's Weyl law \cite[Theorem 0.1]{Muller01} and existing bounds for the sup. norm of the cusp forms, but the author is unaware of any results of sufficient generality to cover non-spherical forms on a non-compact manifold.
On the other hand, this may be proved quite easily by the method given for $d=2$ below.

For $d=3$, the Poincar\'e series \eqref{eq:MainPoincare} just fails to converge absolutely; the contour shifting in the Whittaker function encounters a pole at $(s_1,s_2)=(-1+r,-1-r)$ and there is no way to shift the $r$ contour so that the powers on both $y$ coordinates simultaneously exceed 2.
On the other hand, if we modify the kernel of the Poincar\'e series in \eqref{eq:MainPoincare} to include an extra factor $(y_1 y_2)^u$, the spectral expansion will now also have the lifts of the $d \le 2$ forms, but will converge rapidly and uniformly on compact sets in $\Re(u) \ge 0$ by the usual Fourier-type analysis on the Whittaker functions of the spectral basis.
The Bruhat decomposition also converges to a holomorphic function of $u$ by the usual contour shifting as in section \ref{sect:TechnicalWeylLaw}, and we see that \eqref{eq:PreKuzSpectral} holds for $d=3$ in spite of the conditional convergence of the Poincar\'e series.

The true difficulty arises for $d=2$; in this case, even the sum of Kloosterman sums in the Fourier coefficients (just) fails to converge absolutely, due to a pole at $(s_1,s_2)=(-\frac{1}{2}+r,-\frac{1}{2}-r)$ in the contour shifting.
It will follow from the arithmetic Kuznetsov formula \cite{ArithKuzII} that the naive $d=2$ Kuznetsov formula still holds, with the sum of Kloosterman sums converging in the conditional sense, but first one needs a reasonable upper bound on the Weyl law, so we prove:
\begin{prop}
For $T > 2$,
\[ \sum_{\substack{\varphi\in\mathcal{S}^{2*}_3\\ \abs{r_\varphi\pm iT} \le 1}} \frac{\abs{\rho^*_\varphi(1,1)}^2}{\cosmu^2(r_\varphi)} \ll T^3. \]
\end{prop}
\begin{proof}
We consider the $L^2$-norm of the vector-valued Poincar\'e series
\begin{align*}
	P_F(g) =& \sum_{\gamma\in U(\Z)\backslash\Gamma} F(\gamma g), \\
	F(xyk) :=& f(xy) \Matrix{0&0&1&0&0} \WigDMat{2}(k), \\
	f(xyk) :=& \psi_{1,1}(x) (2\pi y_1)^{1+s_1} (2\pi y_2)^{1+s_2} \exp(-2\pi(y_1+y_2)),
\end{align*}
where $s_1=\wbar{s_2}=\sigma\pm iT$ with $\sigma > 1$ to be chosen later.

Then if $\varphi \in \mathcal{S}_3^{2*}$ has spectral parameters $\mu(r)$, $r\in i\R$, we have
\begin{align*}
	\innerprod{P_F,\varphi} =& 2\pi^3 \sqrt{6} \Lambda^*(r) \rho^*_\varphi(1,1) \wtilde{F}(r), \\
	\wtilde{F}(r) :=& \frac{1}{\Lambda^*(r)} \int_{\Re(s')=0} \Gamma(s_1-s'_1) \Gamma(s_2-s'_2) \Gamma\paren{\tfrac{1}{2}+s'_1-r} \Gamma\paren{\tfrac{1}{2}+s'_2+r} \\
	& \qquad \times B\paren{\tfrac{2+s'_1+2r}{2},\tfrac{s'_2-2r}{2}} \frac{ds'}{(2\pi i)^2}.
\end{align*}
Now shift the $s'$ contours back to $-1+\epsilon$ and suppose $\abs{r\pm iT} \le 1$, so $\wtilde{F} = \wtilde{F}_1+\wtilde{F}_2+\wtilde{F}_3$ with $\wtilde{F}_1$ having the shifted contours, $\wtilde{F}_2$ having the residue at $s_2'=2r$ (and the $s_1'$ contour at, say, $\Re(s_1')=0$), and $\wtilde{F}_3$ having the double residue at $(s_1',s_2')=(-\frac{1}{2}+r,-\frac{1}{2}-r)$.
Then $\wtilde{F}_1(r) \ll T^{-2+\epsilon}$, $\wtilde{F}_2(r)$ decays exponentially in $T$, and $\wtilde{F}(r) \asymp\wtilde{F}_1(r) \asymp T^{-\frac{3}{2}}$.

By general $L^2$ theory, we may complete the orthonormal set of cusp forms $\mathcal{S}_3^{2*}$ to some basis (which won't be a spectral basis) of the whole $L^2$ space and apply positivity to Plancherel's theorem:
\[ \sum_{\abs{r_\varphi\pm iT} \le 1} \frac{\abs{\rho^*_\varphi(1,1)}^2}{\cosmu^2(r_\varphi)} \ll T^3 \innerprod{P_F,P_F}. \]
On the other hand,
\begin{align*}
	\innerprod{P_F,P_F} \le& 4 \sum_{w\in W} \sum_{c_1,c_2\ge 1} S_w(1,1,c) F_w(c), \\
	F_w(c) :=& \int_{Y^+} \int_{\wbar{U}_w(\R)} \abs{f(cwxt)} dx \, (2\pi t_1)^{1+\sigma} (2\pi t_2)^{1+\sigma} \exp(-2\pi(t_1+t_2)) dt.
\end{align*}
We conjugate $xt=tu$ and for each $w\in W$, we choose some $\beta=\beta(w)\in\R^3$ so $\exp(-y_1-y_2) \ll p_{-\beta(w)}(y)$ and
\[ f(cwtu) p_{\rho+\sigma\rho}(t) \abs{\frac{dx}{du}} \ll p_{\rho+\sigma\rho^w+\sigma\rho}(t) p_{\rho+\sigma\rho}(cwu) p_{-\beta}(cwtu) \ll p_{3\rho}(cwu) t_1^{h_1} t_2^{h_2} \]
with $h_1,h_2 > 2$ also depending on $w$.
This is sufficient for convergence of the $x$ and $t$ integrals and the $c$ sum so that $\innerprod{P_F,P_F} \ll 1$.

Taking $\sigma=4$, a valid choice of $\beta(w)$ is
\[ \beta(I) = \beta(w_2) = \beta(w_3) = 0, \qquad \beta(w_4)=\beta(w_5)=\beta(w_l)=2\rho. \]
\end{proof}

Note that the above method only relied on the evaluation of the Mellin transform of a single entry of the vector-valued Whittaker function; in fact, this particular entry was not only the simplest in form, but also the easiest to evaluate \cite[sect. 7.2]{HWII}, and it is reasonable to expect this type of behavior to exist on a wide variety of groups.
So we expect that the above method would generalize nicely to give (weak, but simple) bounds on arithmetically-weighted Weyl laws in those situations, as well.
This is also interesting from the $L$-function standpoint as it gives a lower bound on the adjoint-square $L$-function at 1 (aka. the residue of the Rankin-Selberg $L$-function).

\bibliographystyle{amsplain}

\bibliography{HigherWeight}

\end{document}